\documentclass[12pt]{article}
\usepackage{amssymb}
\usepackage{amsmath}
\usepackage{theorem}
\usepackage{latexsym}
\usepackage{psfig}
\usepackage{graphicx}
\theoremheaderfont{\bf} 
\theorembodyfont{\sl}
\newtheorem{theo}{Theorem}[section]
{\theorembodyfont{\rm} \newtheorem{defi}[theo]{Definition}}
{\theorembodyfont{\rm} \newtheorem{exa}[theo]{Example}}
{\theorembodyfont{\rm} \newtheorem{rem}[theo]{Remark}}
{\theorembodyfont{\rm} }
{\theorembodyfont{\rm} }

\newenvironment{proof}{{\sc Proof:}}{\mbox{}\hfill$\Box$\par}

\newcommand{\Section}[1]{\section{#1}\setcounter{equation}{0}}

\newcommand{\junk}[1]{}

\newcounter{abc}

\title{A Block Decomposition Algorithm for Computing Rook~Polynomials
\date{June 4, 2004}
\author{
  Abigail Mitchell\\
  {\small Department of Mathematics}\vspace{-2mm}\\
  {\small University of Notre Dame}\vspace{-2mm}\\
  {\small Notre Dame, IN 46556-5683, USA}\vspace{-2mm}\\
  {\small {\em e-mail:} amitche3@nd.edu} 
   }}

\begin{document}
\maketitle
\begin{abstract}
  Rook polynomials are a powerful tool in the theory of restricted
  permutations.  It is known that the rook polynomial of any board
  can be computed recursively, using a cell decomposition technique
  of Riordan. \cite{ri58}

  In this paper, we give a new decomposition theorem,
  which yields a more efficient algorithm for computing the rook
  polynomial.  We show that, in the worst case, this block decomposition
  algorithm is equivalent to Riordan's method.\bigskip

\end{abstract}

\clearpage
\Section{Introduction}

Rook polynomials provide a method of enumerating permutations with
restricted position.  Their study was begun in 1946 by Kaplansky
and Riordan~\cite{ka46r}, with applications to card-matching problems.
Riordan's 1958 book~\cite{ri58} is considered the first systematic
analysis, and remains a classic treatment of the subject.  A series
of papers by Goldman et al. \cite{go75jw,go76jrw,go78jw,go77jw,go76jw} in the 1970s expanded the field by
applying more advanced combinatorial methods.  

More recently, Haglund \cite{ha96a,ha98os} has made
investigation into various connections of rook polynomials to other
parts of mathematics: hypergeometric series, enumeration of matrices
over finite fields, and group representation theory.  Rook polynomials
are also closely related to matching theory, chromatic theory and
various other graph-theoretic topics (for example, \cite{go78jw,fa91w,ch95g}.)  Applications to quantum mechanics have
recently begun to surface, through the use of rook polynomials with
Weyl algebras \cite{va04}.  In combinatorics proper, rook polynomials have been related to various permutation statistics \cite{bu04}, and the inverse problem has been solved for Ferrers boards \cite{mi04a}.

It has long been known that the rook polynomial of any board can be
computed recursively.  In this paper, we present a new block decomposition
algorithm for computing rook polynomials.  This is a fundamental
generalization of Riordan's cell decomposition \cite{ri58}, and in
general improves on its efficiency.

\section{Preliminaries}\label{sec:prelim}

Let $B$ be a generalized chessboard -- a set of cells arranged in rows
and columns.  The rook is a chess piece which attacks on rows and
columns; by a rook placement we shall mean a non-attacking placement
of $k$ indistinguishable rooks on the board $B$.

This intuitive definition can be formalized in several ways.
A board $B$ may be regarded as a subset of $[1,2,\ldots ,m]\times [1,2,\ldots ,n]$
for some $m,n \in \mathbb{N} $.
  A rook placement on $B$ then corresponds to a choice of $k$ elements
of $B$, ${(x_1,y_1),(x_2,y_2),\ldots ,(x_k,y_k)}$ such that ${x_i}={x_j}$ or
${y_i}={y_j}$ implies $i=j$.  

Alternatively, let $D\subseteq [1,2,\ldots ,m]$ be an arbitrary subset.  Then a rook placement is an injective function
 $f:D \rightarrow [1,2,\ldots ,n]$
 such that, for any $i\in D$, $(i,f(i))\in B$.

$B$ may also be seen as an element of ${M_{m,n}}[\mathbb{F}_2]$, an
$m \times n$ matrix with binary entries.  We write in this case $B=(b_{i,j})$.  In this interpretation, a
placement of $k$ rooks on $B$ corresponds to a choice of $k$
independent 1's in $B$.

A third formalization is to regard $B$ as a bipartite graph on vertex
sets $[1,2,...,m]$ and $[1,2,...,n]$, where the graph contains edge
$(i,j)$ iff $(i,j)$ is a cell in the corresponding chessboard.  A
rook placement then corresponds to a partial matching on the graph.

\begin{defi}
The rook polynomial of a board $B$ is the ordinary generating function
\[
R(B;x)=\sum_{k=0}^{\infty}r_k(B) x^k,
\]
where the coefficient $r_k(B)$ is the number of placements of $k$ rooks
on $B$.  Where its omission will cause no confusion, we will assume
the variable to be $x$, and will simply write $R(B)$.
\end{defi}

Rook polynomials were first investigated by Kaplansky and Riordan
\cite{ka46r} in the 1940's.  The following preliminary results and
observations, which we shall state here without proof, are theirs;
proofs can be found in \cite{ri58}.

\begin{defi}
Two boards $A$, $B$ are said to be disjoint if no cell of $B$ is in
the same row or column as any cell of $A$.
\end{defi}

\begin{rem}\label{dbt}
In this case, it is clear
that the number of ways of placing $k$ rooks on a board consisting of the disjoint union $A\amalg B$
is
\[
r_k(A \amalg B) =\sum_{i=0}^{k} r_i(A) r_{k-i}(B),
\]
where $r_i(A)$ and $r_i(B)$ are coefficients in the rook polynomials of $A$
and $B$, respectively.  Then the rook polynomial of $A\amalg B$ is
\[
R(A\amalg B;x)=R(A;x)R(B;x).
\]
\end{rem}

\begin{theo}\label{rbt}
The rook polynomial of a rectangular $m \times n$ board, denoted $R_{m,n}(x)$
or simply $R_{m,n}$, is given by
\[
R_{m,n}(x)=\sum_{k=0}^{\min(m,n)} {m\choose k} {n\choose k} k! \, x^k.
\]
\end{theo}

\begin{figure}
\begin{center}
  \includegraphics[totalheight=0.5in]{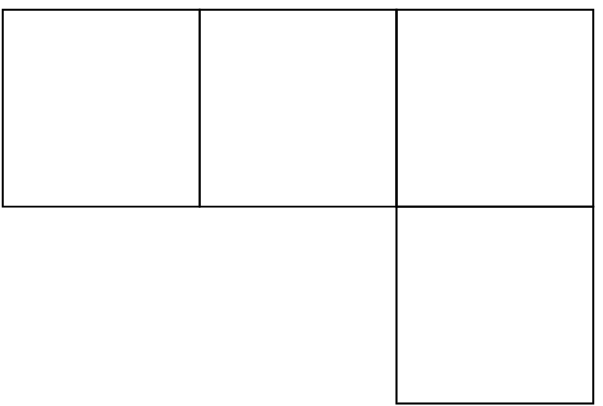}
  \hspace{1.0in}
  \includegraphics[totalheight=0.5in]{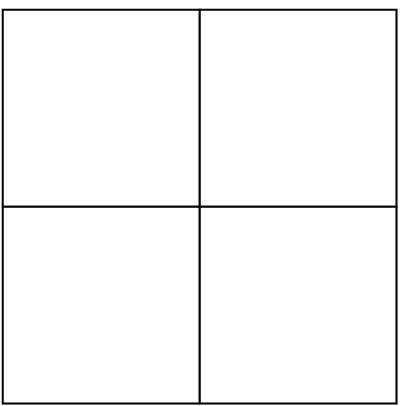}
\caption{Two rook equivalent boards, with rook polynomial $1+4x+2x^2$.} \label{equivalent-boards}
\end{center}
\end{figure}

Two boards $A$,$B$ are said to be rook equivalent if $R(A)=R(B)$.  A
sufficient condition for rook equivalence is that $B$ can be obtained
from $A$ by permutation of rows and columns.  (This condition is not
necessary, as demonstrated by the two boards in Figure~\ref{equivalent-boards}.)

Riordan \cite{ri58} also gives an algorithm which allows the rook
polynomial of any board to be computed recursively.

\begin{theo}[Cell Decomposition Theorem]\label{cdt}
Let $B$ be a board, and $b_{i,j}$ be a cell of $B$.  Let $B_e$ denote the board obtained from $B$ by deleting cell $b_{i,j}$, and let $B_i$ denote the board obtained from $B$ by deleting row $i$ and column $j$.  Then
\[
R(B)=R(B_e)+xR(B_i)
\]
\end{theo}

This decomposition theorem provides a recursive method for calculating the rook polynomial of any board.  In the next
section, we present a generalization of this algorithm, which allows
more efficient computation.

\section{Block decomposition algorithm}\label{sec:main}

In this section, we shall view a board formally as an $m \times n$
matrix over $\mathbb{F}_2$.

\begin{defi}
A subboard $S$ of a board $B$ is a set of cells of $B$, inheriting the row and column relationships of $B$.  Intuitively, $S$ is, up to permutation of rows and columns, a submatrix of $B$.  More formally, $S$ is represented by a matrix $(s_{i,j})$ over $\mathbb{F}_2$, together with injective mappings $\phi_1$ and $\phi_2$ from the row and column indices of $S$ to those of $B$, having the property that $s_{i,j}=b_{\phi_1(i),\phi_2(j)}$.
\end{defi}

\begin{defi}
A subboard $S$ of $B$ is said to cover a row $i$ of $B$ (respectively a column $j$) if $i$ (resp. $j$) is contained in the image of $\phi_1$ (resp. $\phi_2$.)
\end{defi}

\begin{defi}
A block $S$ in $B$ is a subboard of $B$ satisfying the following conditions:
\begin{enumerate}
\item For any rows $i,i'$ covered by $S$, and any column $j$ not covered by $S$, $b_{i,j}=b_{i',j}$.
\item For any row $i$ not covered by $S$, and any columns $j,j'$ covered by $S$, $b_{i,j}=b_{i,j'}$.
\end{enumerate}
\end{defi}  

Intuitively, an $s \times t$ block $S$ is a subboard of $B$ formed at the
intersection of a set of $s$ rows with a set of $t$
columns, all of which are identical except perhaps on $S$.  (Note that every cell of $B$ is trivially a 1$\times$1 block.)

\begin{defi}
Let $B$ be an $m \times n$ board, and $S$ an $s \times t$
block in $B$. For $0\leq j\leq \min(s,t)$, let $B_{S,j}$ denote
the board obtained from $B$ by deleting
\begin{enumerate}
\addtolength{\itemsep}{-0.5\baselineskip}
\item $j$ of the $s$ rows covered by $S$
\item $j$ of the $t$ columns covered by $S$
\item all the cells of $S$
\end{enumerate}
$B_{S,j}$ is called the $j^{th}$ inclusion board of $B$ relative to $S$.  Note that this is well-defined, since the rows and columns involved are identical except on $S$, which is also deleted.
\end{defi}

The following theorem is a generalization of Riordan's cell decomposition, Theorem~\ref{cdt}.

\begin{theo}[Block Decomposition Theorem]\label{bdt}
Let $B$ be an $m \times n$ board, and $S$ an $s \times t$ block in $B$.
Let $r_k(S)$ be the coefficient of $x^k$ in the rook polynomial of $S$.
Let $B_{S,j}$ be the $j^{\mbox{th}}$ inclusion board of $B$ relative to
$S$, for $1\leq j \leq \min(s,t)$.
Then
\[
R(B)=\sum_{j=0}^{\min(s,t)}r_j(S)\, x^j R(B_{S,j})
\]
\end{theo}

\begin{proof}
Consider all the $r_k(B)$ possible placements of $k$ rooks on $B$.  These
may be partitioned into $k+1$ classes, according to the number of rooks
placed on $S$.  Rook placements which include $j$ rooks on $S$ are
exactly those enumerated by the $j^{\mbox{th}}$ inclusion board
$B_{S,j}$.  Hence, the number of placements of $k$ rooks on $B$ is
\[
r_k(B) = \sum_{j=0}^{k} r_j(S) r_{k-j}(B_{S,j})].
\]
Therefore,
\[
R(B)=\sum_{j=0}^{\infty} r_j(S) x^j R(B_{S,j}).
\]
\end{proof}

\section{Remarks on efficiency}\label{efficiency}

As noted above, any cell $b_{i,j}$ of $B$ is a $1\times 1$ block.  By
Theorem~\ref{rbt}, its rook polynomial is $1+x$.  Therefore,
Theorem~\ref{bdt} gives
\[
R(B)=R(B_{S,0})+xR(B_{S,1}).
\]

$B_{S,0}$ and $B_{S,1}$ are exactly the exclusion and inclusion boards
$B_e$ and $B_i$ defined by Riordan (as given in Theorem~\ref{cdt}),
so that in the case of a $1\times 1$ block, this decomposition agrees
with Riordan's cell decomposition.

Since $1\times 1$ is the smallest block possible, and since such a block always exists, this block decomposition algorithm will always be at least as fast as Riordan's cell decomposition.  Efficiency of either algorithm is greatly affected by the choice of cell/block by which to decompose, as shown in the next few examples.

\begin{exa}
Consider the right-triangular board shown below.  The decompositions given below both result from $2\times2$ blocks: there are four such blocks to choose from.  As shown below, a wise choice will immediately reduce the number of subsequent iterations needed by providing more nearly-rectangular boards, and therefore larger blocks available in the next iterations.

\begin{center}
  \includegraphics[totalheight=0.72in]{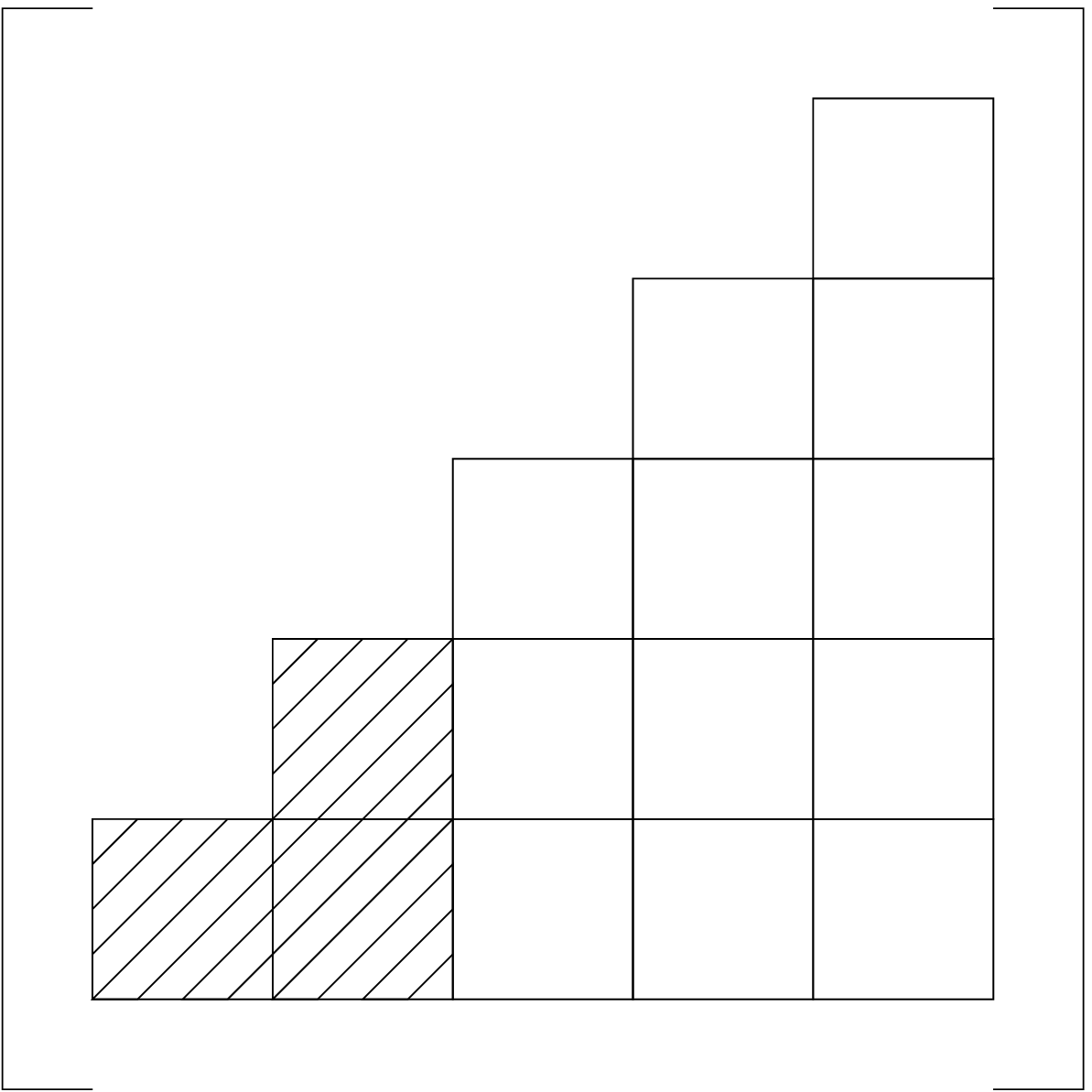}
  \raisebox{0.31in} =
  \includegraphics[totalheight=0.72in]{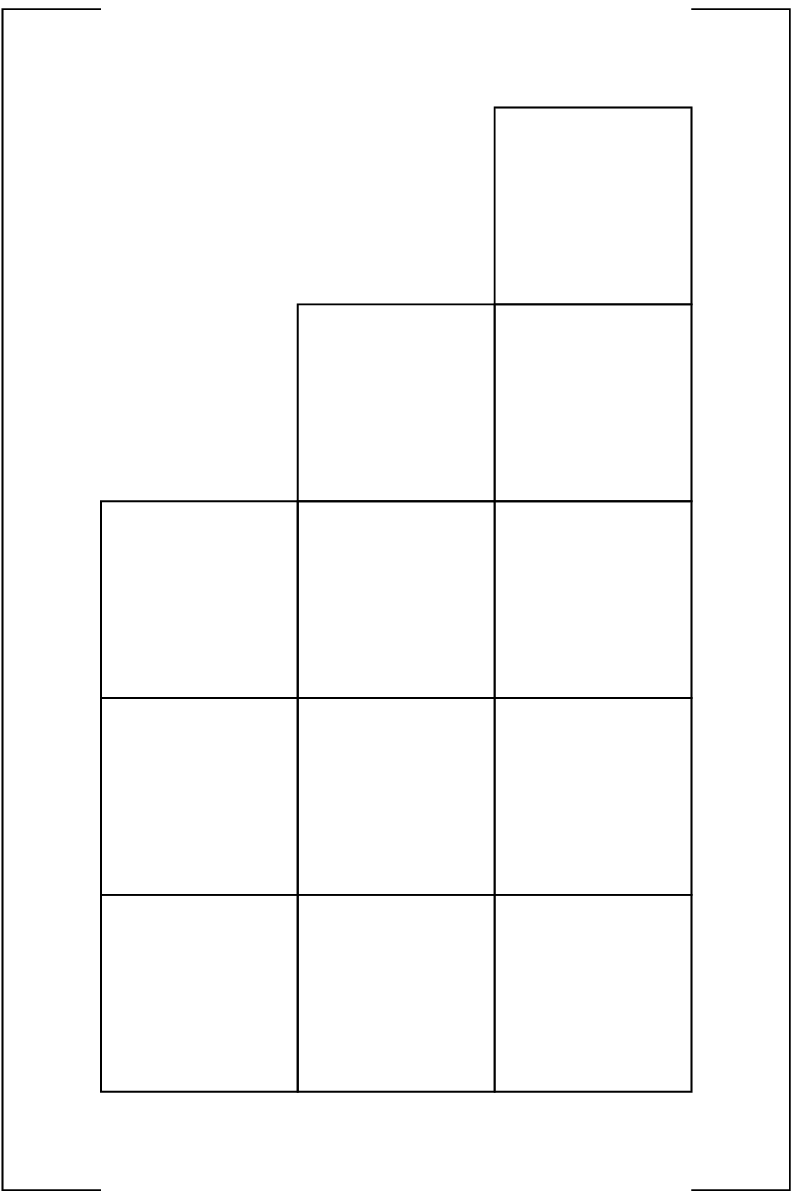}
  \raisebox{0.31in} {$+\;3x$}~\raisebox{0.06in}{\includegraphics[totalheight=0.60in]{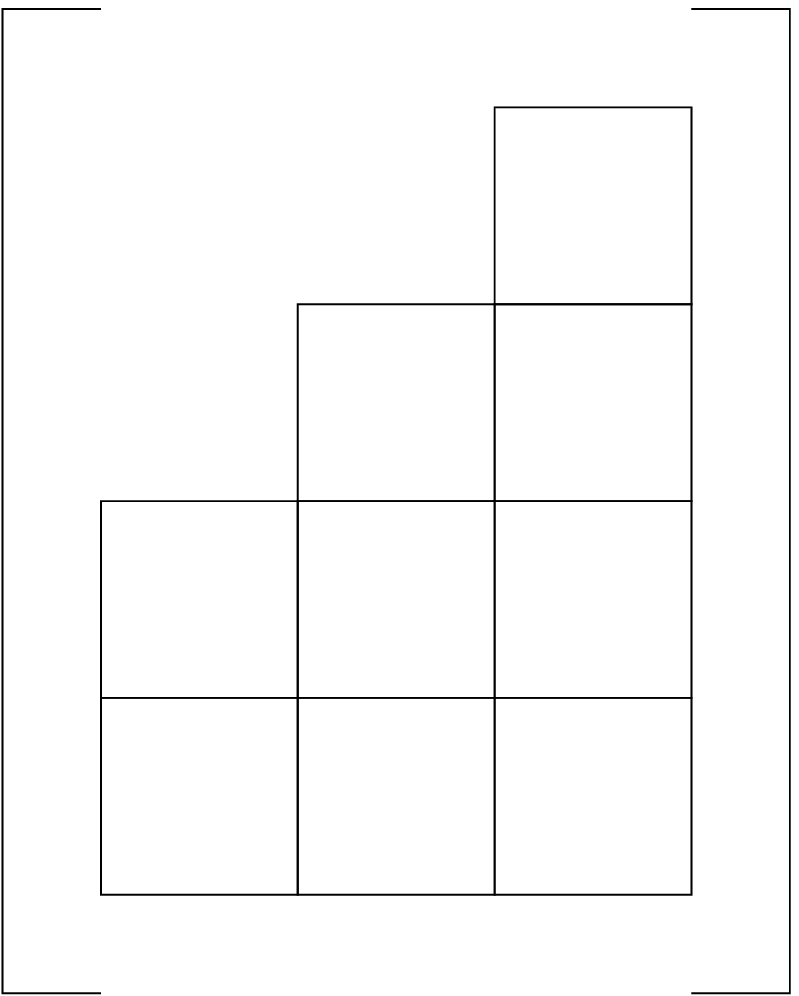}}
  \raisebox{0.31in} {$+\;x^2$}~\raisebox{0.12in}{\includegraphics[totalheight=0.48in]{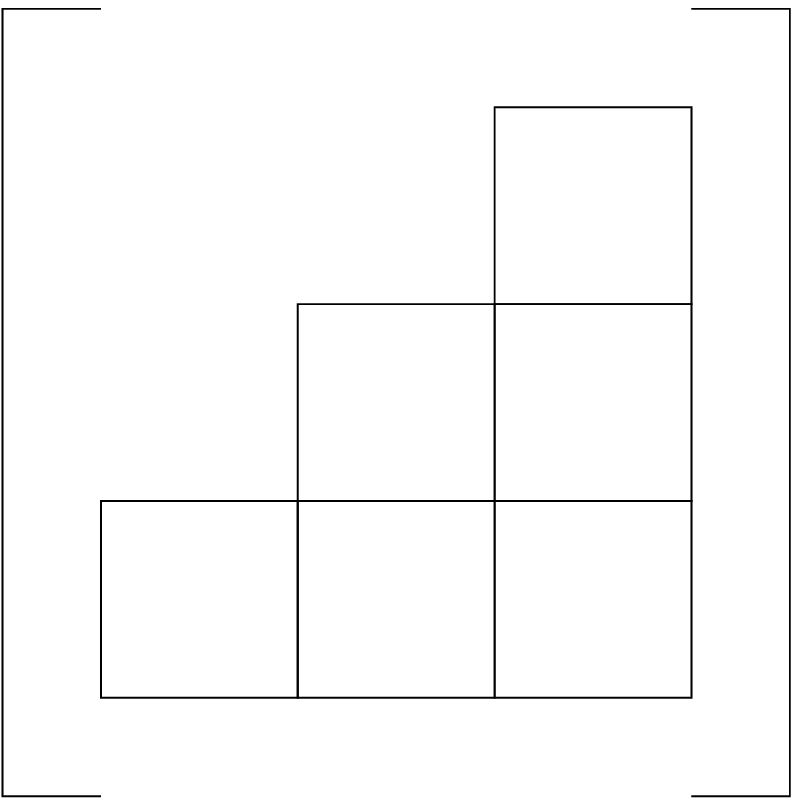}}
  \raisebox{0.31in} {\hspace{0.5in} (good choice)}
  
  \vspace{0.2in}
  \includegraphics[totalheight=0.72in]{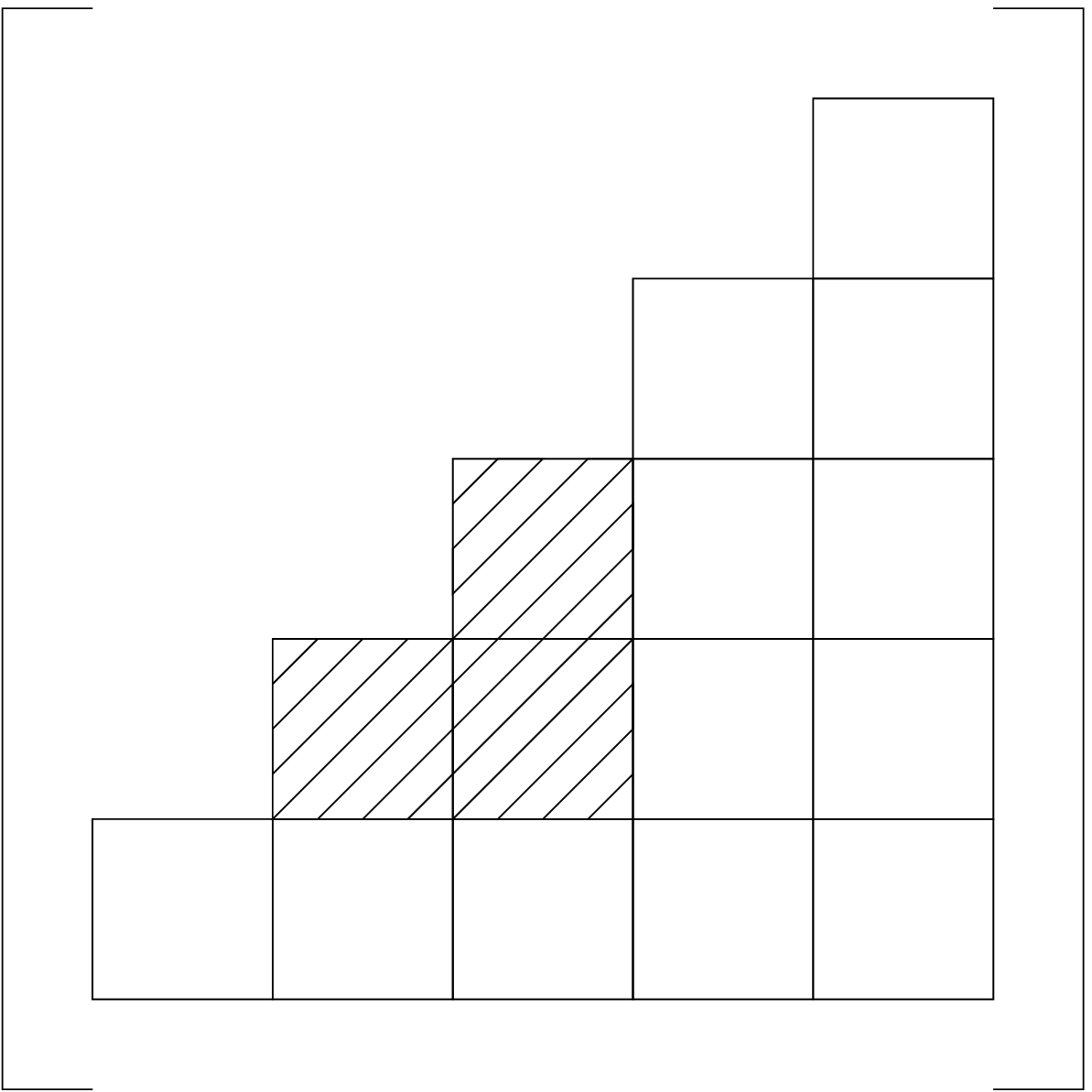}
  \raisebox{0.31in} =
  \includegraphics[totalheight=0.72in]{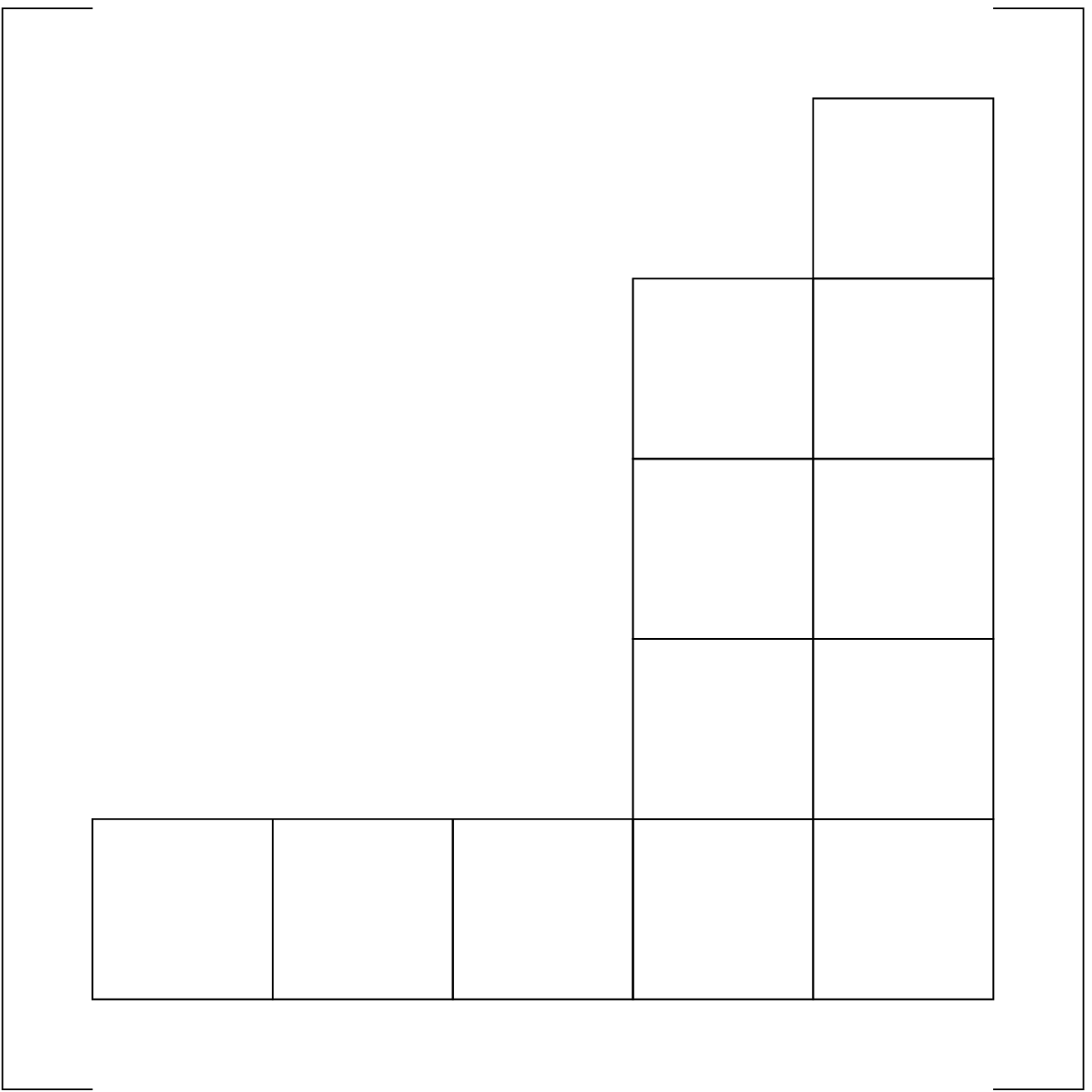}
  \raisebox{0.31in} {$+\;3x$}~\raisebox{0.06in}{\includegraphics[totalheight=0.60in]{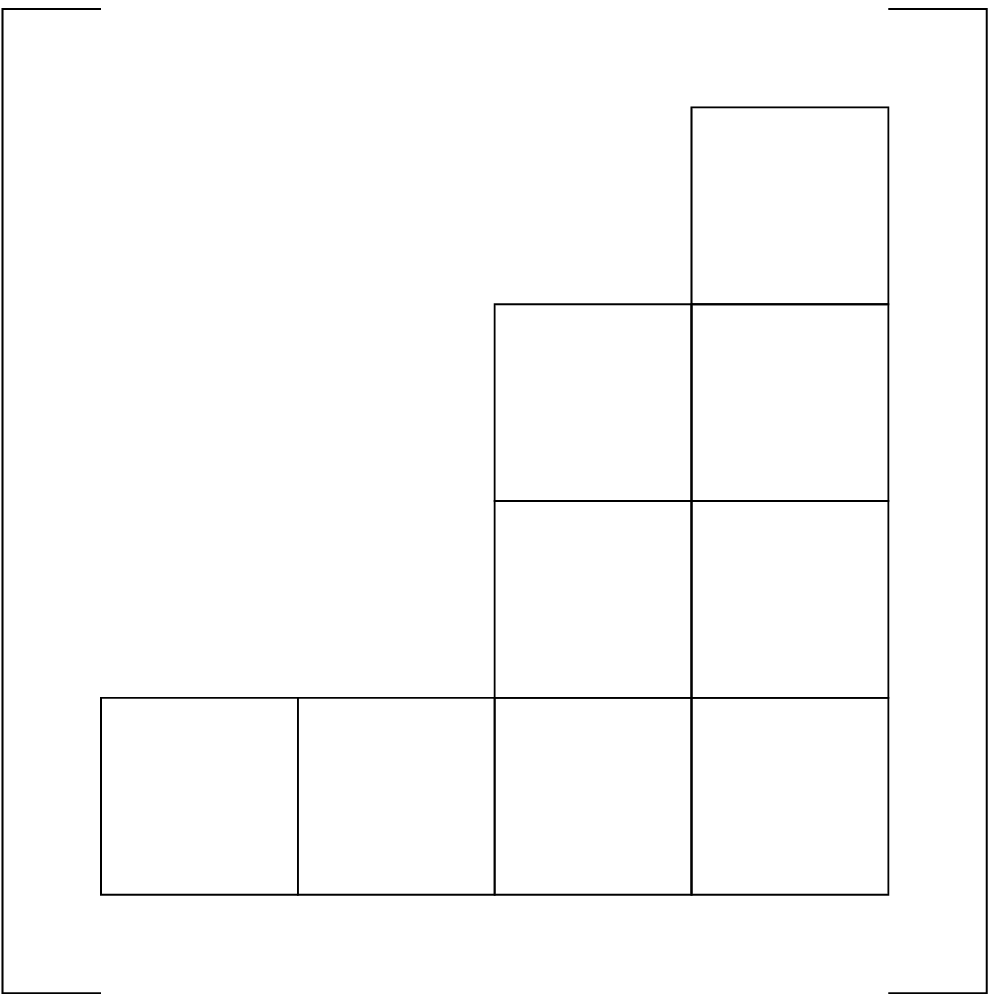}}
  \raisebox{0.31in} {$+\;x^2$}~\raisebox{0.12in}{\includegraphics[totalheight=0.48in]{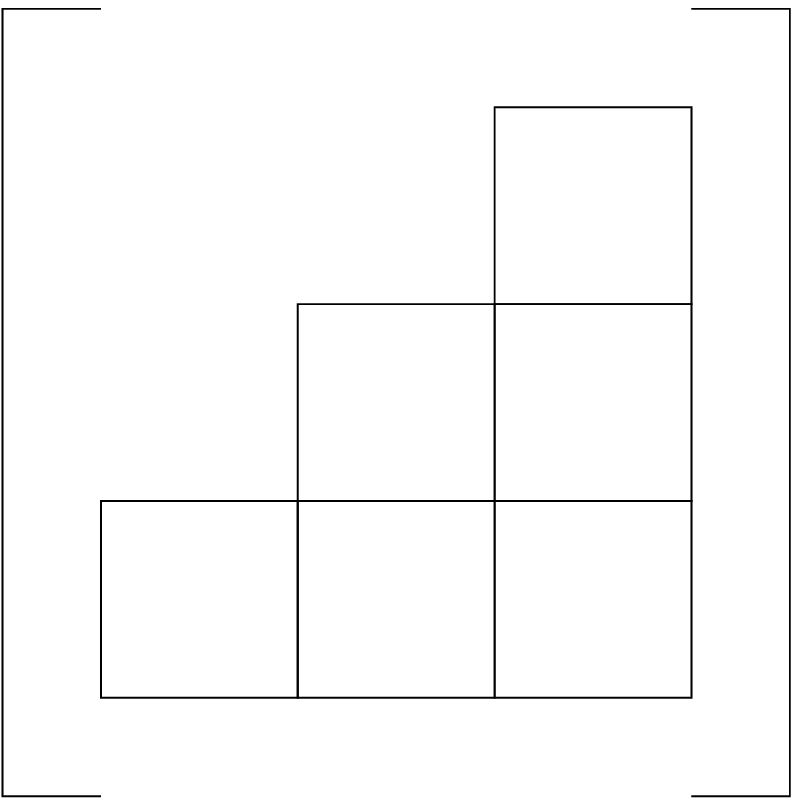}}
  \raisebox{0.31in} {\hspace{0.5in} (poor choice)}
\label{choices}
\end{center}

\end{exa}

Making good block choices becomes even more important for boards where large blocks exist, or boards which are nearly disjoint.

\begin{exa} Consider the following board, which has a $3\times 4$ block composed of the shaded cells shown.
\vspace{0.2in}

\begin{center}
  \includegraphics[totalheight=0.72in]{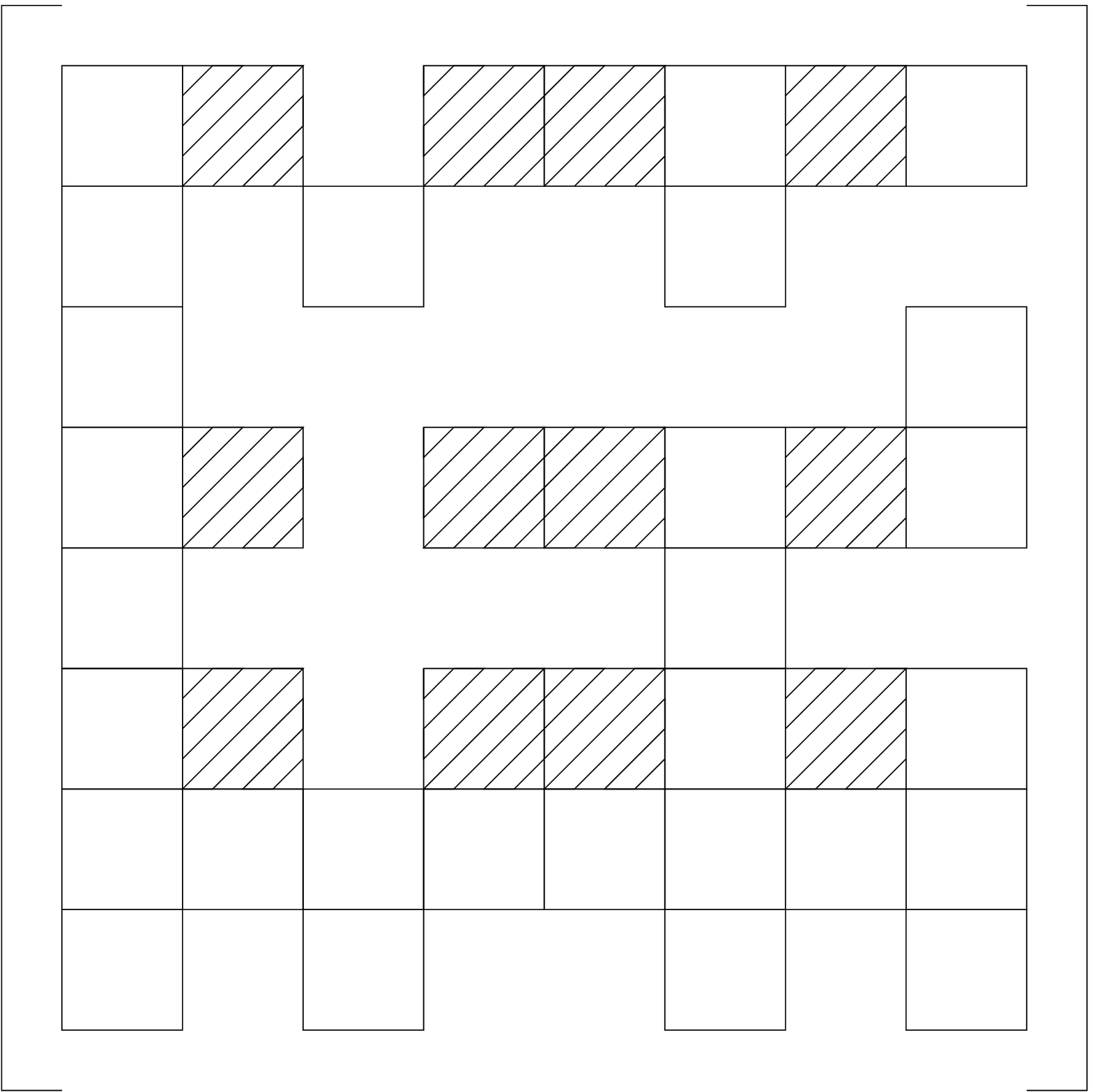}
  \raisebox{0.31in} =
  \includegraphics[totalheight=0.72in]{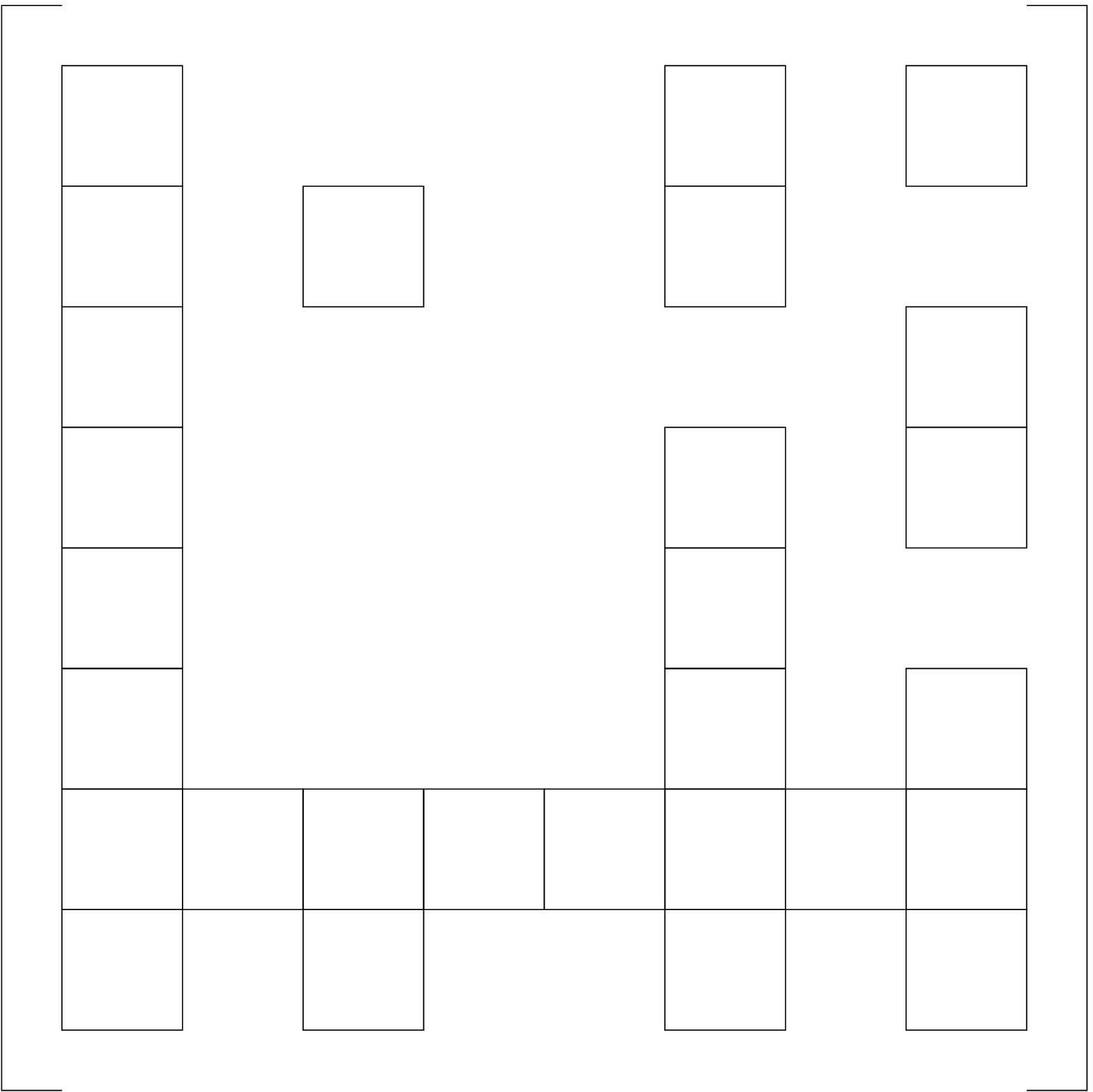}
  \raisebox{0.31in} {$+\;12x$}~\includegraphics[totalheight=0.72in]{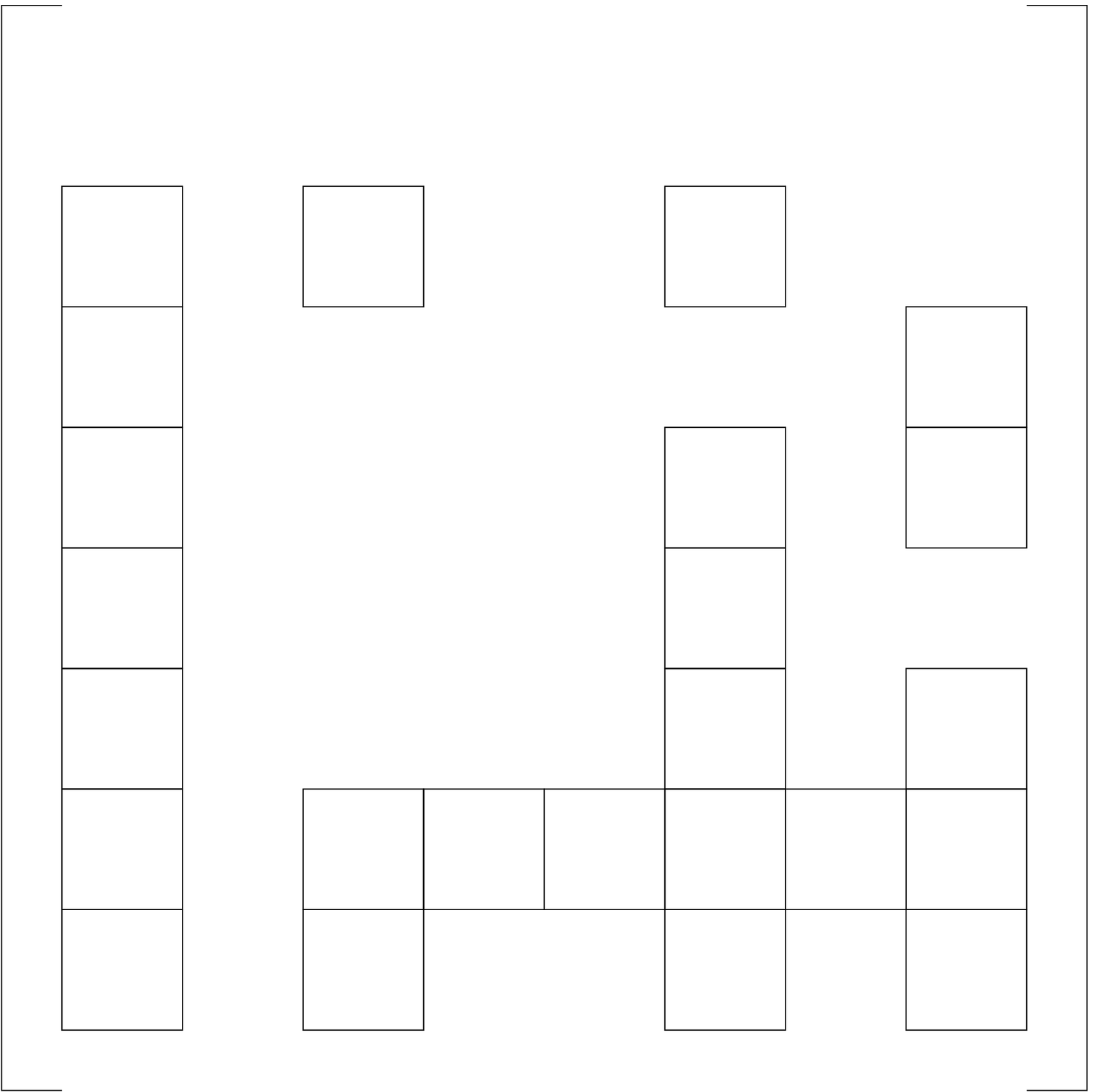}
  \raisebox{0.31in} {$+\;36x^2$}~\includegraphics[totalheight=0.72in]{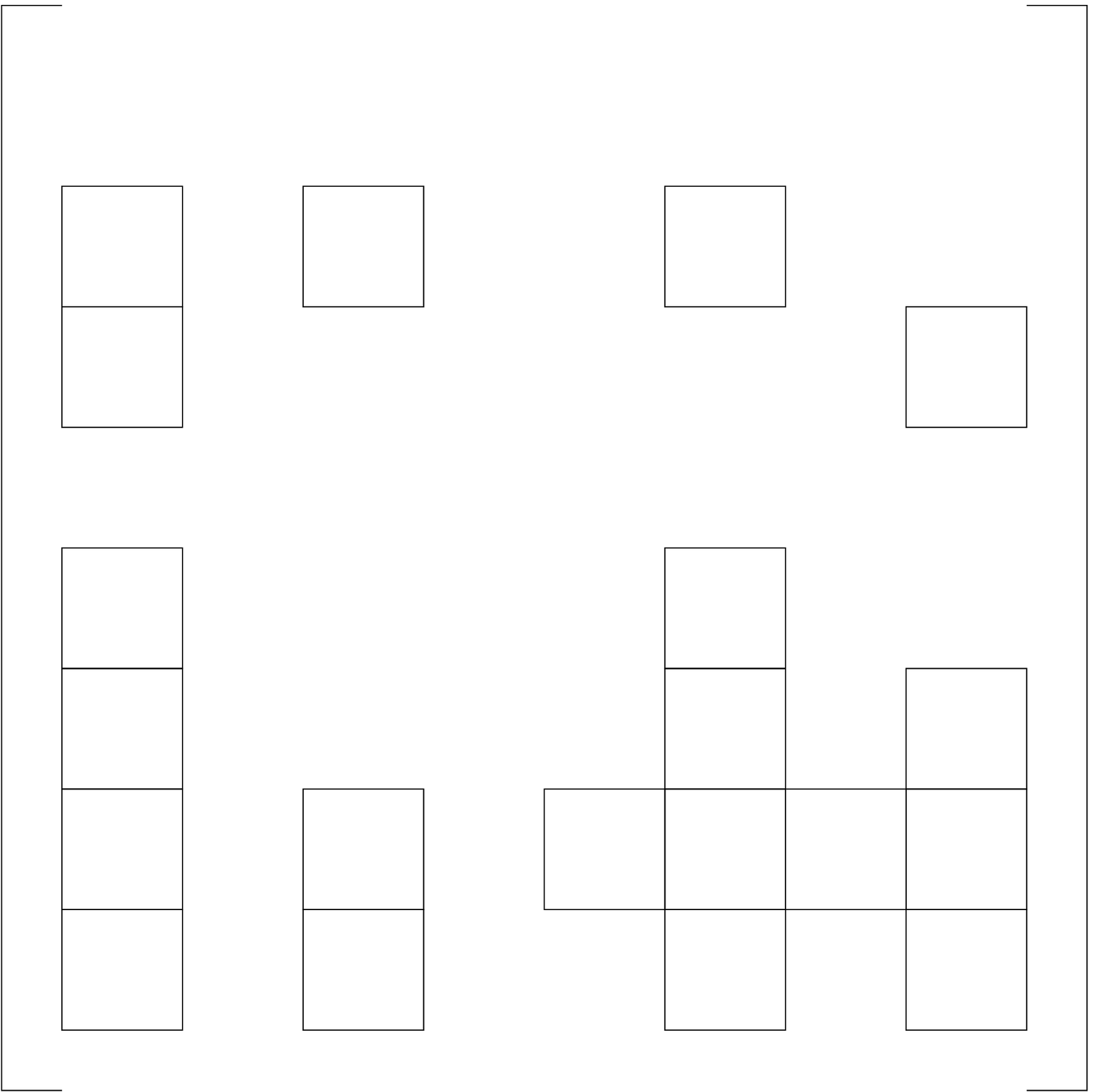}
  \raisebox{0.31in} {$+\;24x^3$}~\includegraphics[totalheight=0.72in]{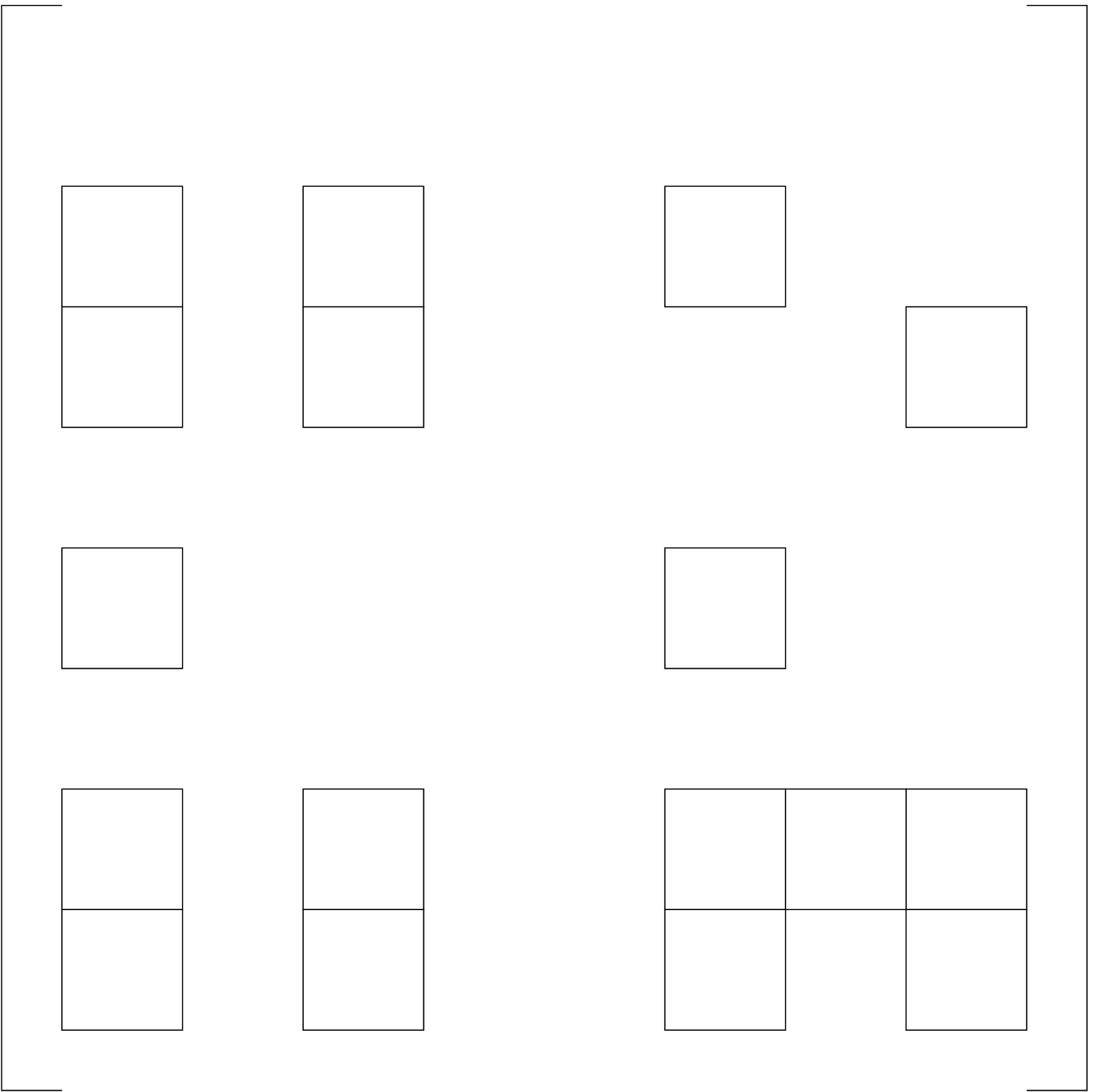}  

\end{center}
After permuting rows and columns, equivalently  
\begin{center}
  \vspace{0.2in}
  \includegraphics[totalheight=0.72in]{bigblock_001.eps}
  \raisebox{0.31in} =
  \includegraphics[totalheight=0.72in]{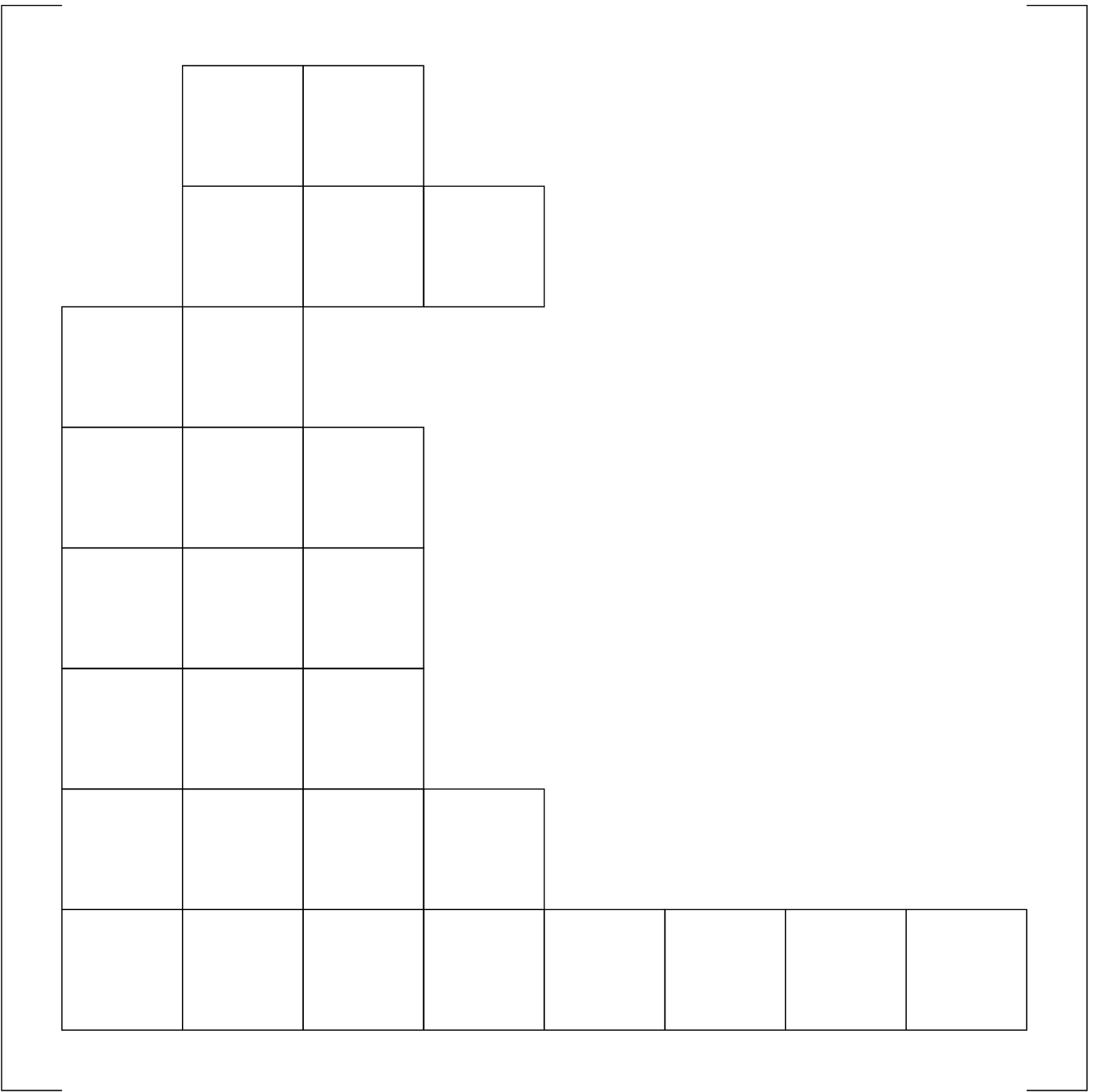}
  \raisebox{0.31in} {$+\;12x$}~\raisebox{0.04in}{\includegraphics[totalheight=0.64in]{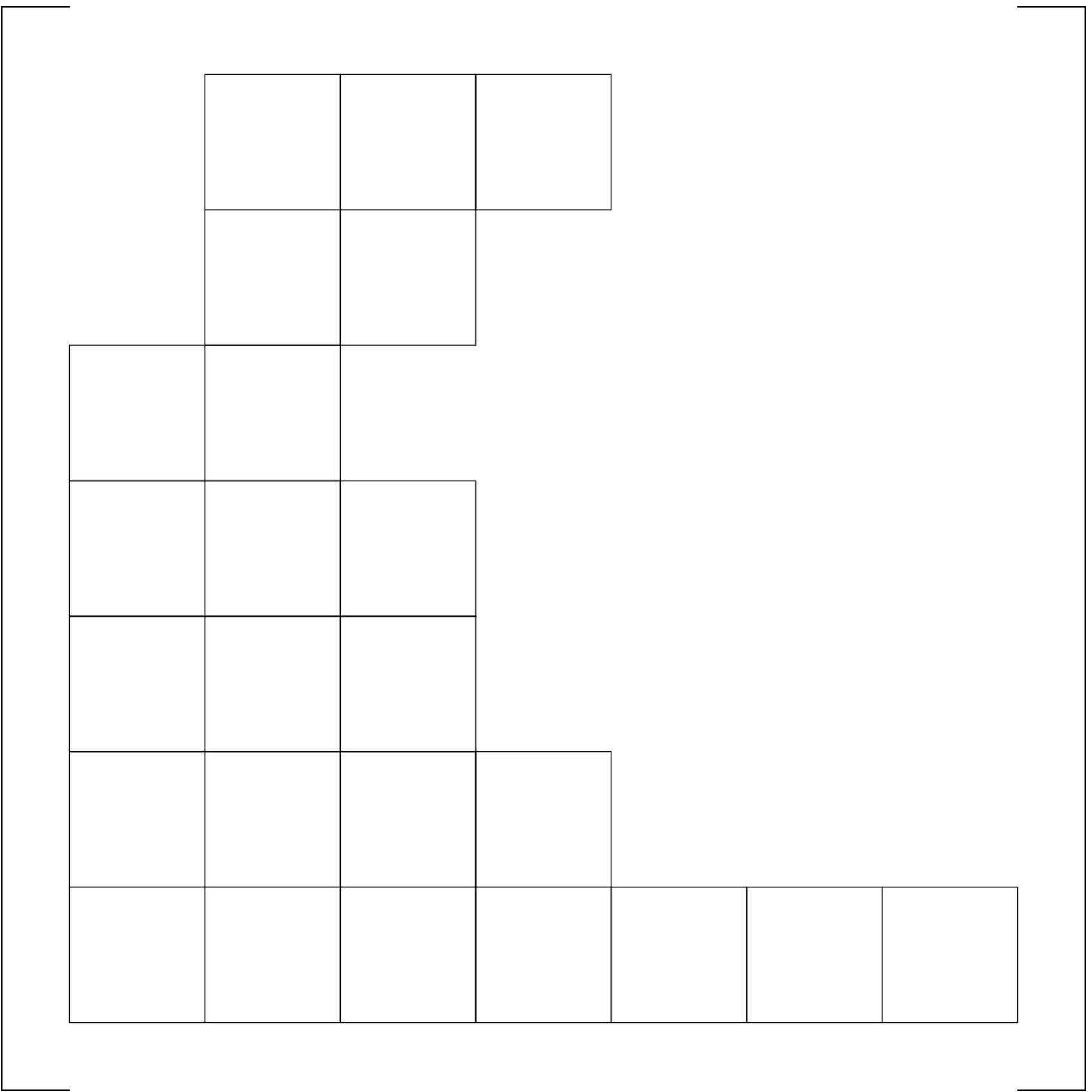}}
  \raisebox{0.31in} {$+\;36x^2$}~\raisebox{0.08in}{\includegraphics[totalheight=0.56in]{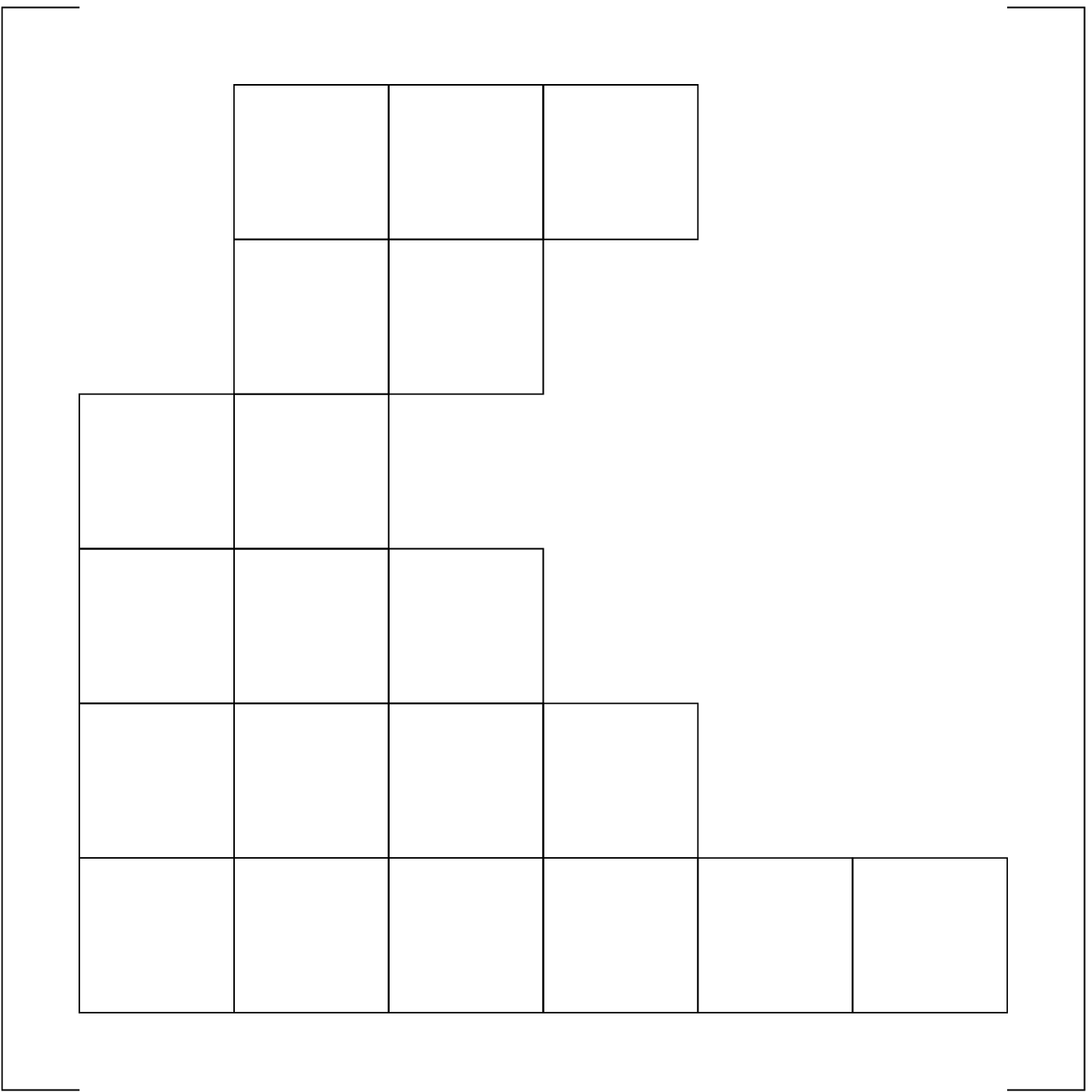}}
  \raisebox{0.31in} {$+\;24x^3$}~\raisebox{0.12in}{\includegraphics[totalheight=0.48in]{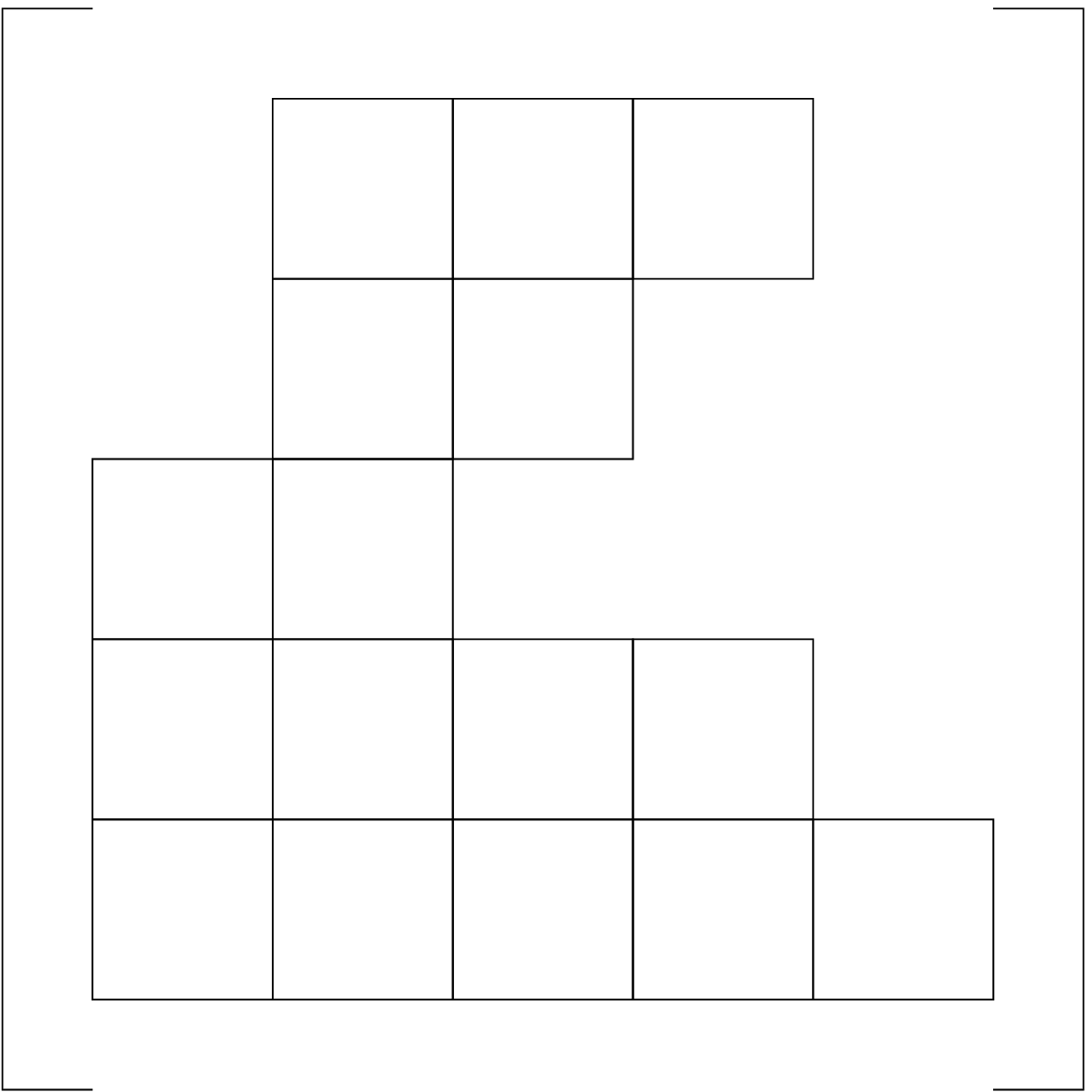}}

\label{bigblock}
\end{center}

Note that we have decomposed 12 cells with one iteration, expressing the desired rook polynomial in terms of four simpler boards.  The cell decomposition method of Theorem\ref{cdt} would have required 12 iterations to achieve the same thing, and the resulting expression would have involved $2^{12}$ boards!

\end{exa}

From the preceding example, it is tempting to conjecture that the optimal solution is to use the largest block possible.  However, this is not always the best strategy, as illustrated by the next example.

\begin{exa} Consider the following board, which is "nearly" disjoint.
\begin{center}
\raisebox{0.45in} {$B\; :=\;$}
\includegraphics[totalheight=1.0in]{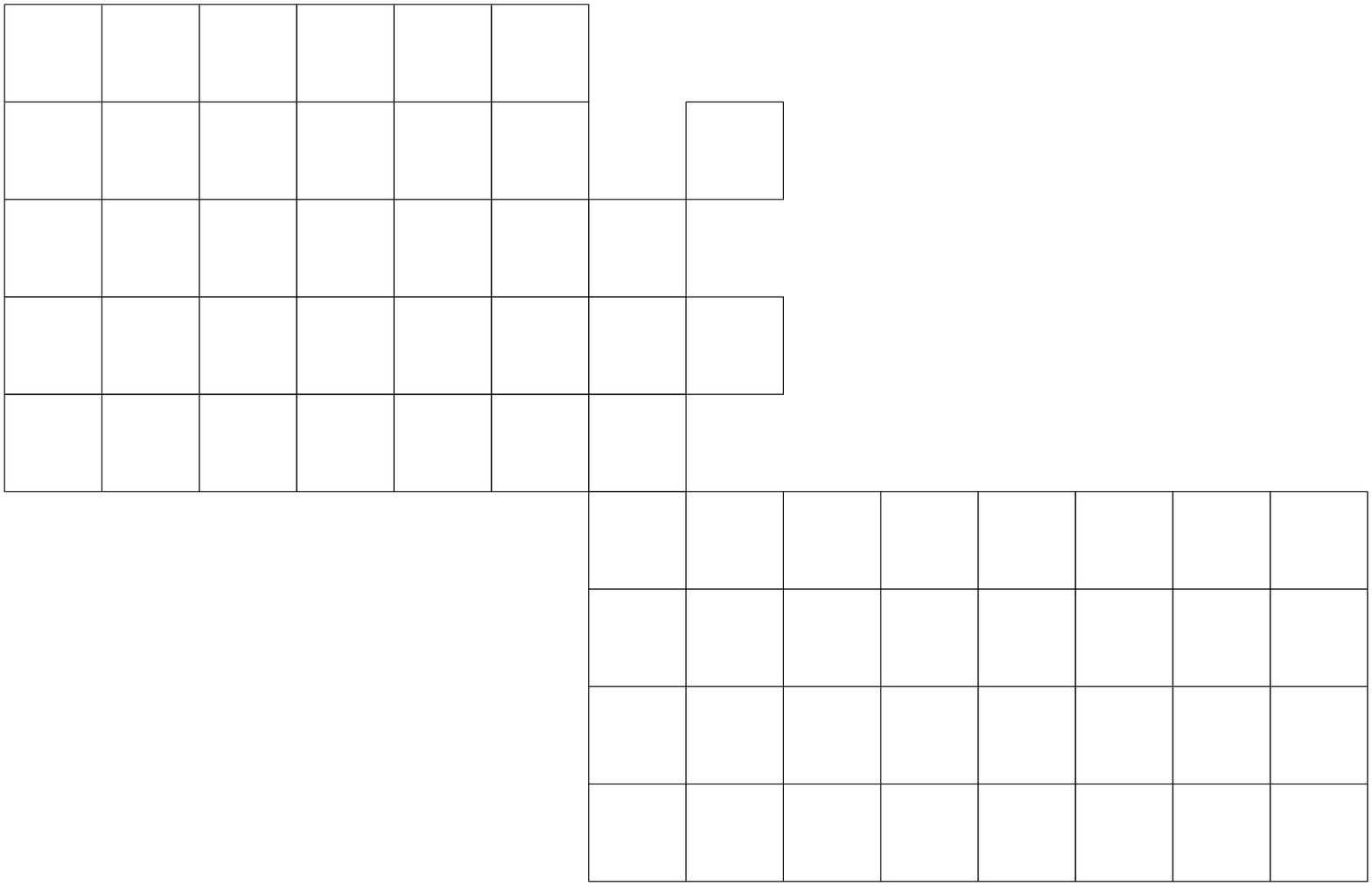}
\end{center}

There are five maximal blocks here: a $1\times 6$ and a $4\times 6$ in the upper left, a $4\times 2$ and a $4\times 6$ in the lower right, and the $4\times 2$ block in the upper right.  The largest is the $4\times 6$; decomposing by this gives

\vspace{0.2in}

\begin{center}
  \includegraphics[totalheight=1.0in]{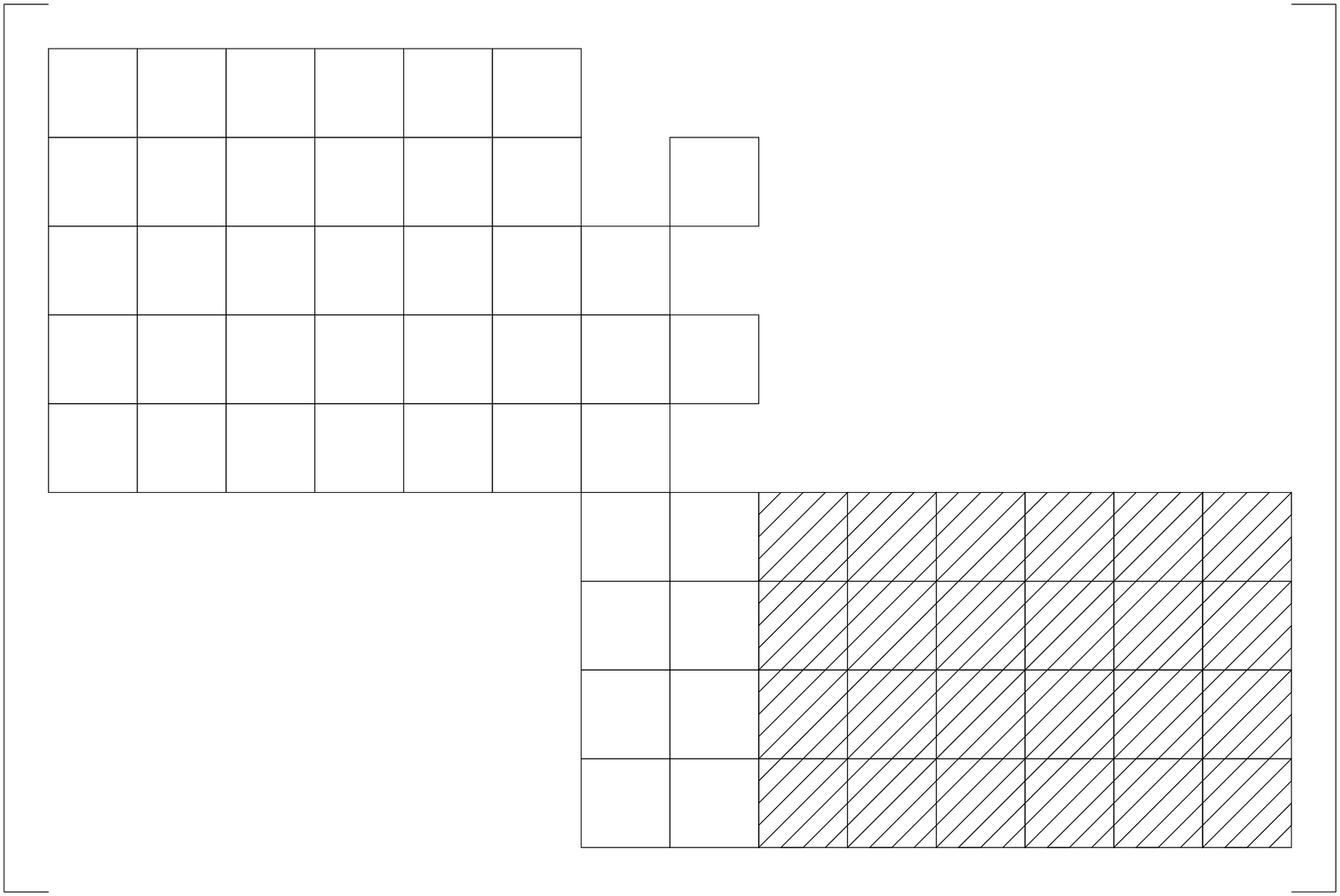}
  \raisebox{0.45in} =
  \includegraphics[totalheight=1.0in]{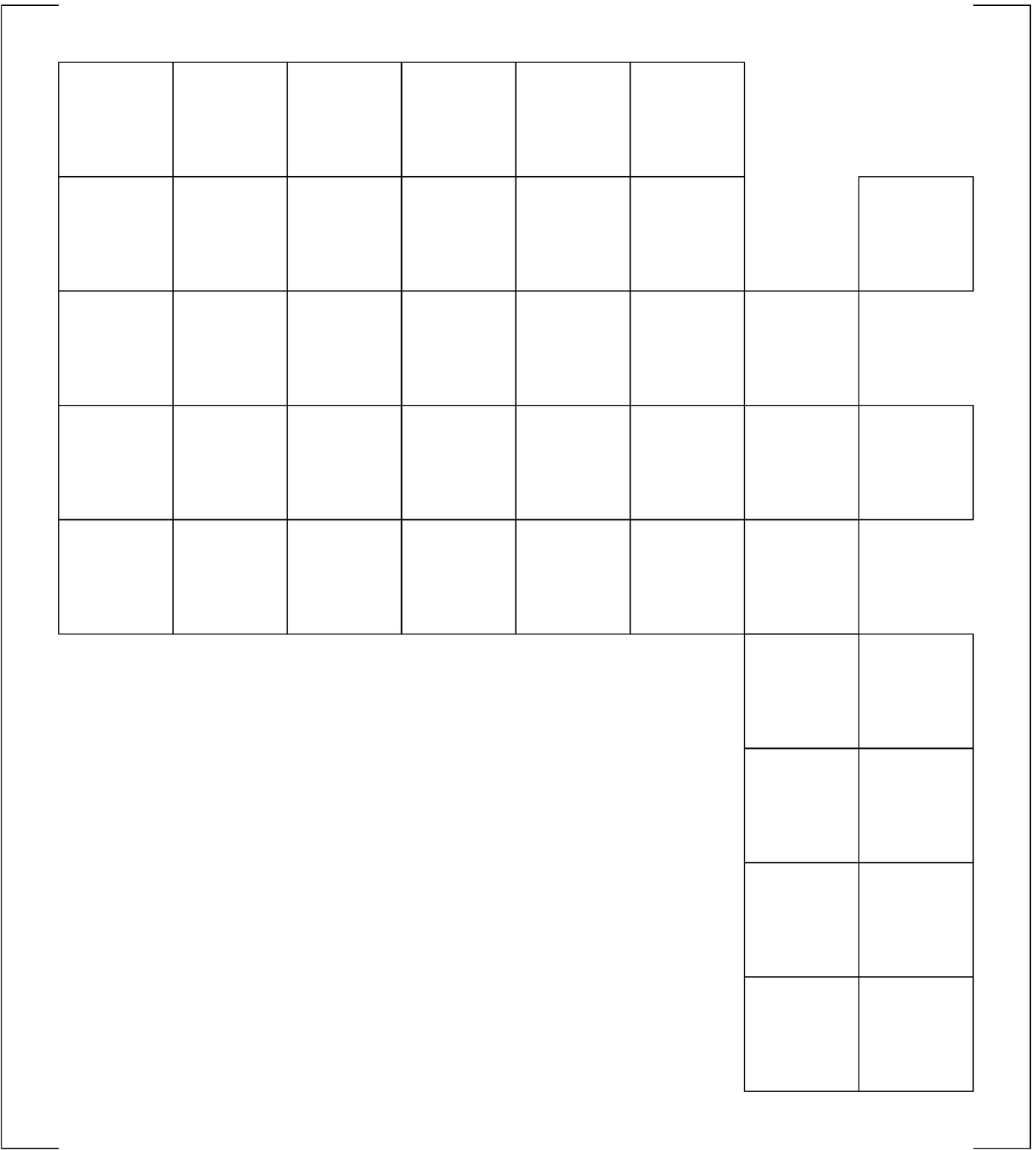}
  \raisebox{0.45in} {$+\;24x$}~\raisebox{0.05in}{\includegraphics[totalheight=0.9in]{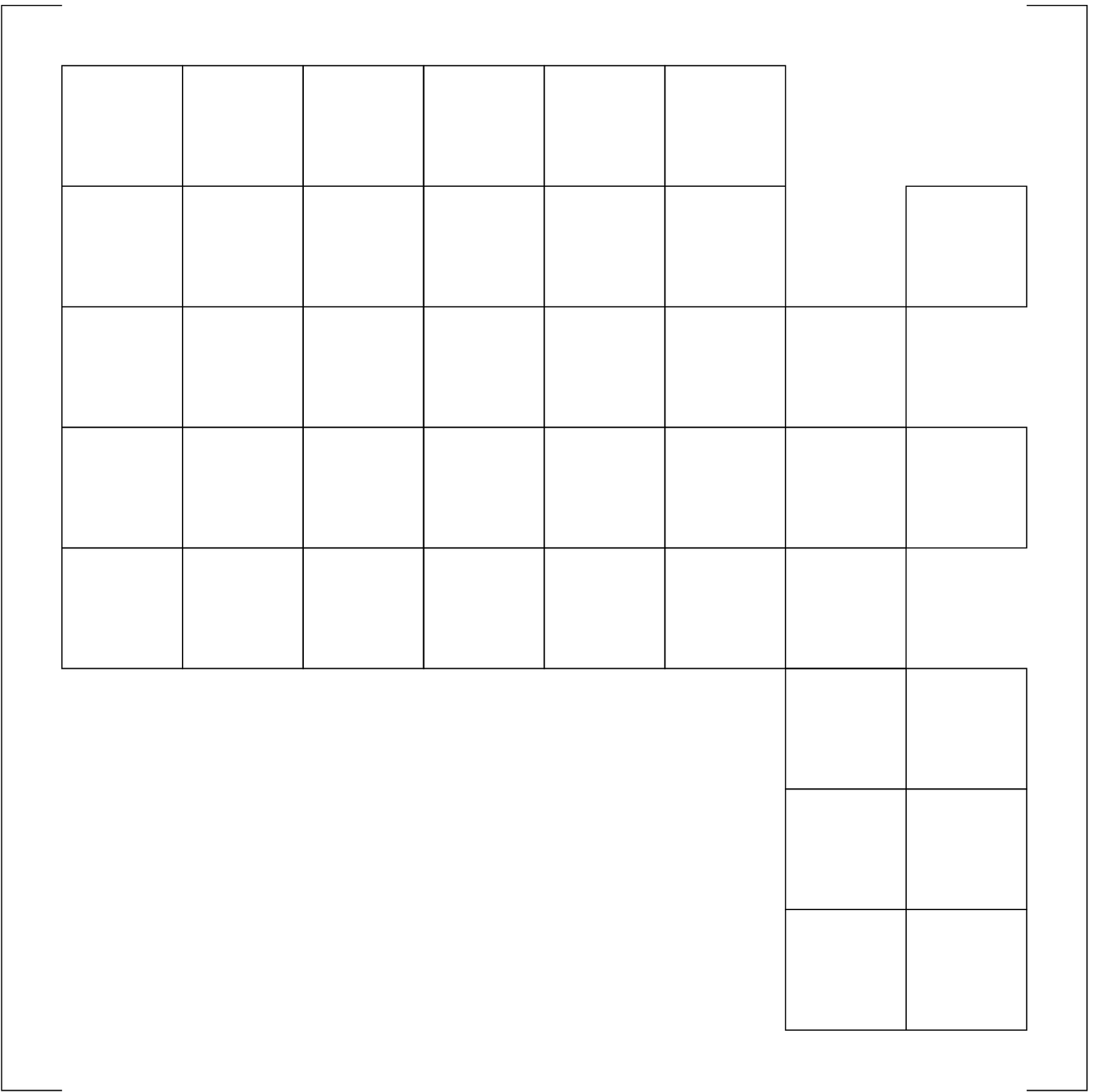}}
  \raisebox{0.45in} {$+\;180x^2$}~\raisebox{0.10in}{\includegraphics[totalheight=0.8in]{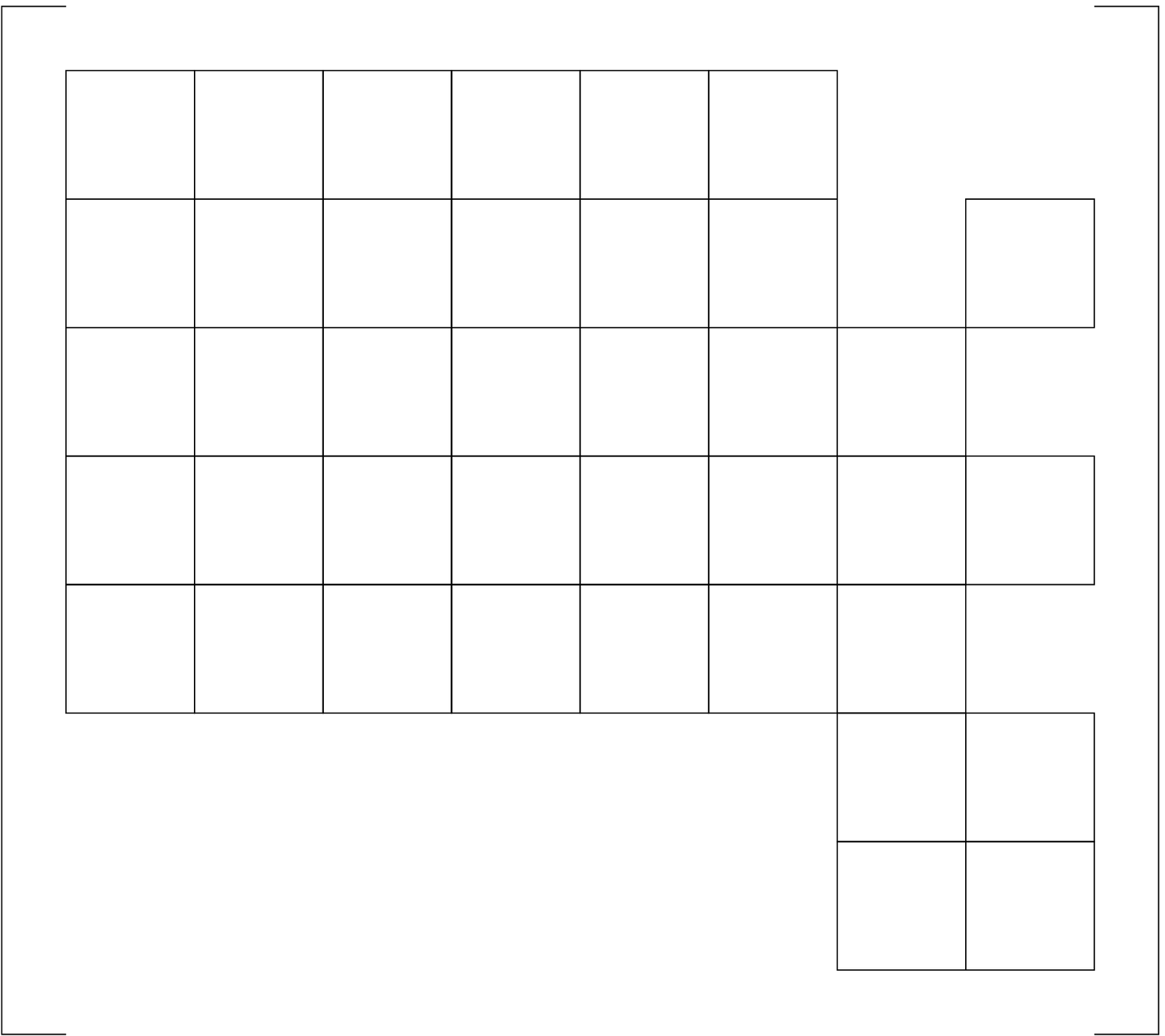}}
  \raisebox{0.45in} {$+\;480x^3$}~\raisebox{0.15in}{\includegraphics[totalheight=0.7in]{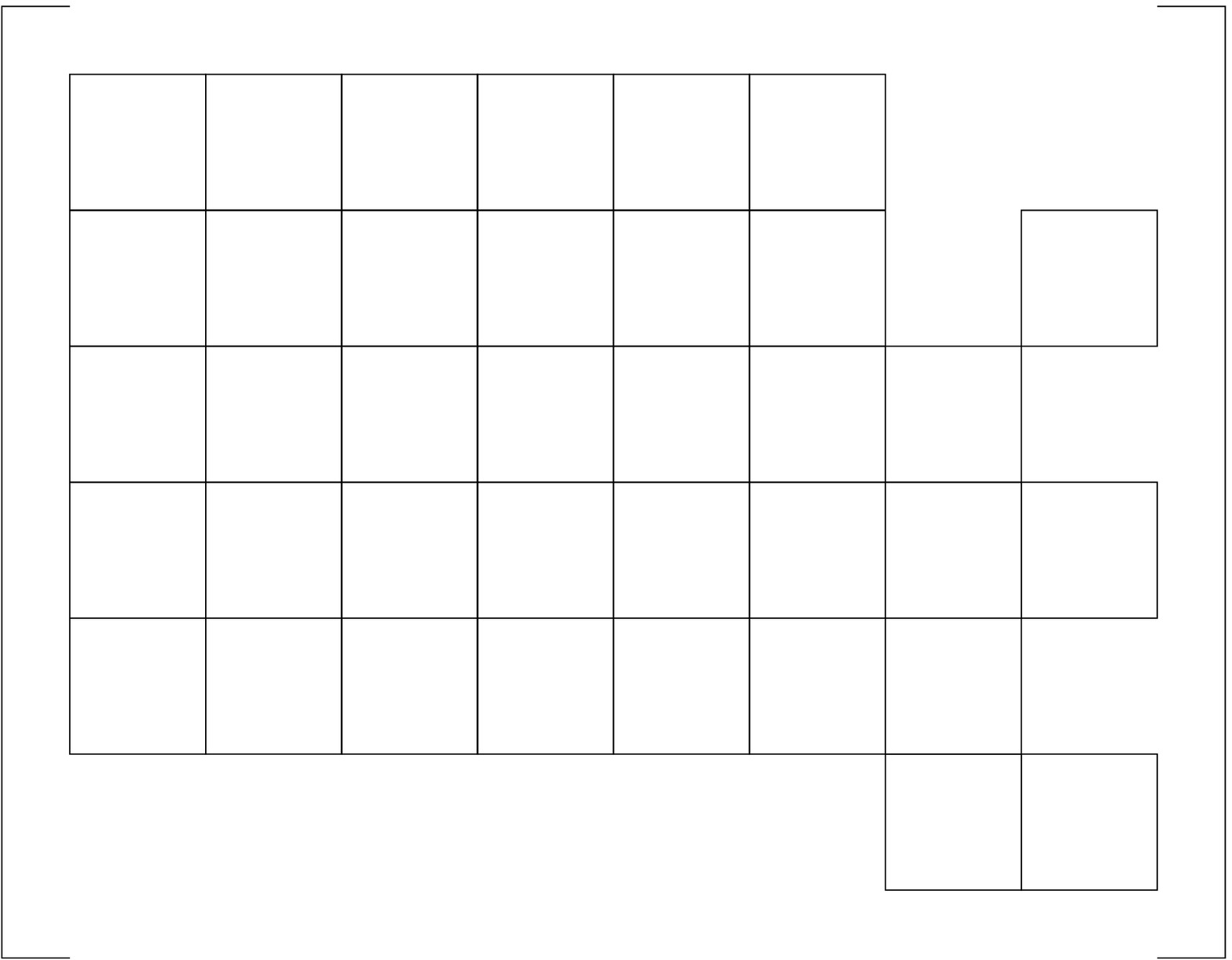}} 
  \raisebox{0.45in} {$+\;360x^4$}~\raisebox{0.20in}{\includegraphics[totalheight=0.6in]{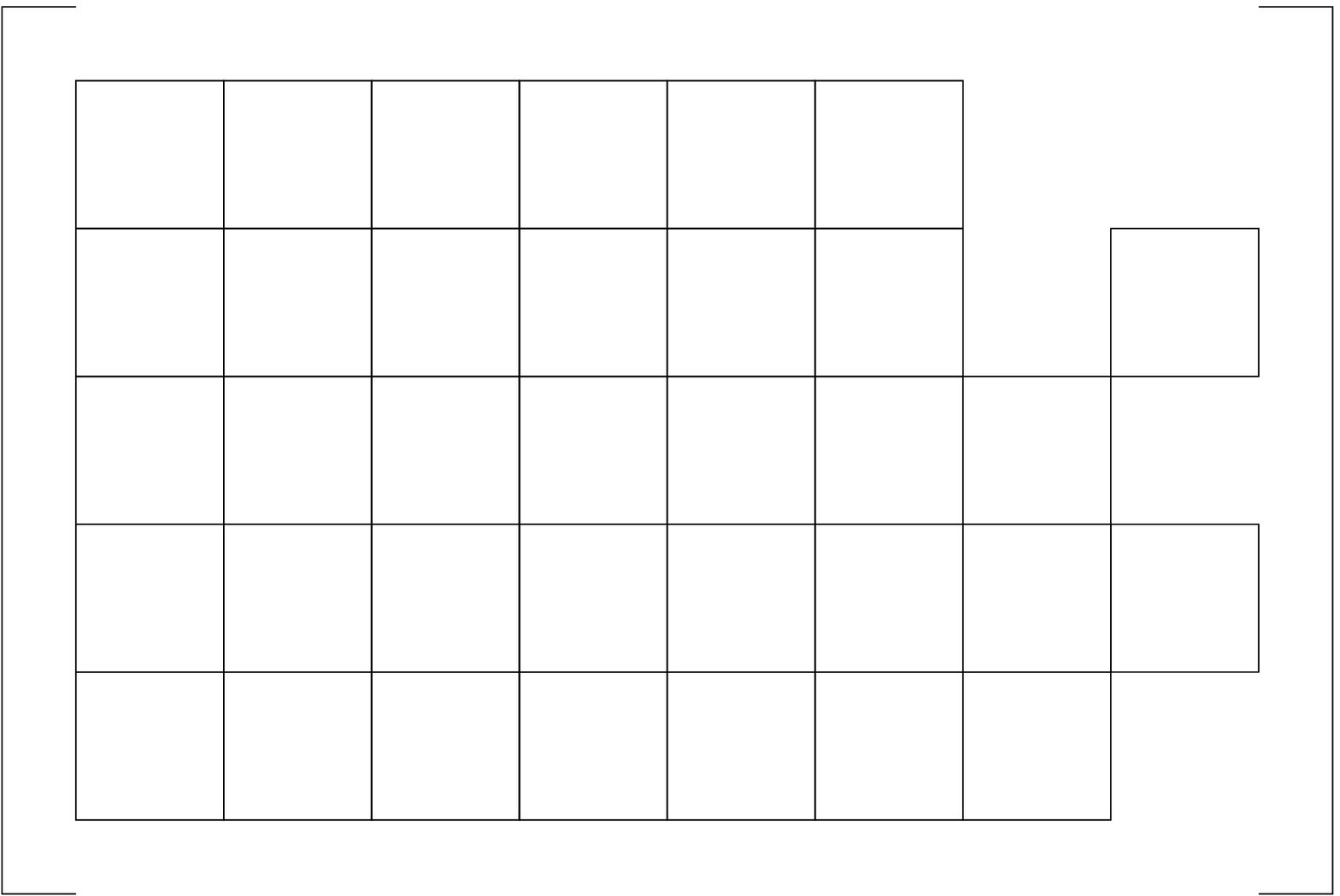}}  
\end{center}

At least one more decomposition iteration will be required on each of these boards.  On the contrary, if we decompose by the $4\times 2$ block,

\begin{center}
  \includegraphics[totalheight=1.0in]{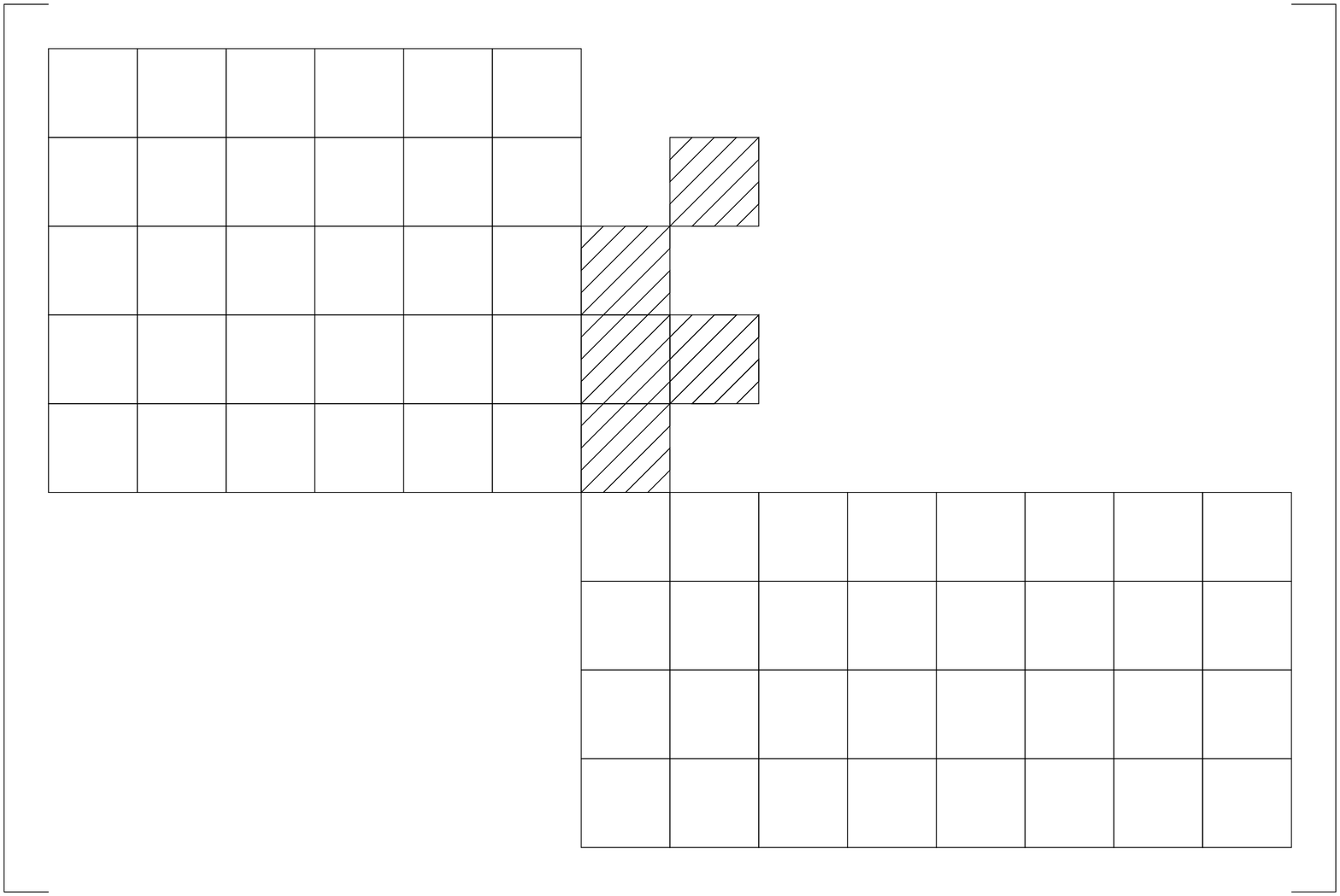}
  \raisebox{0.45in} =
  \includegraphics[totalheight=1.0in]{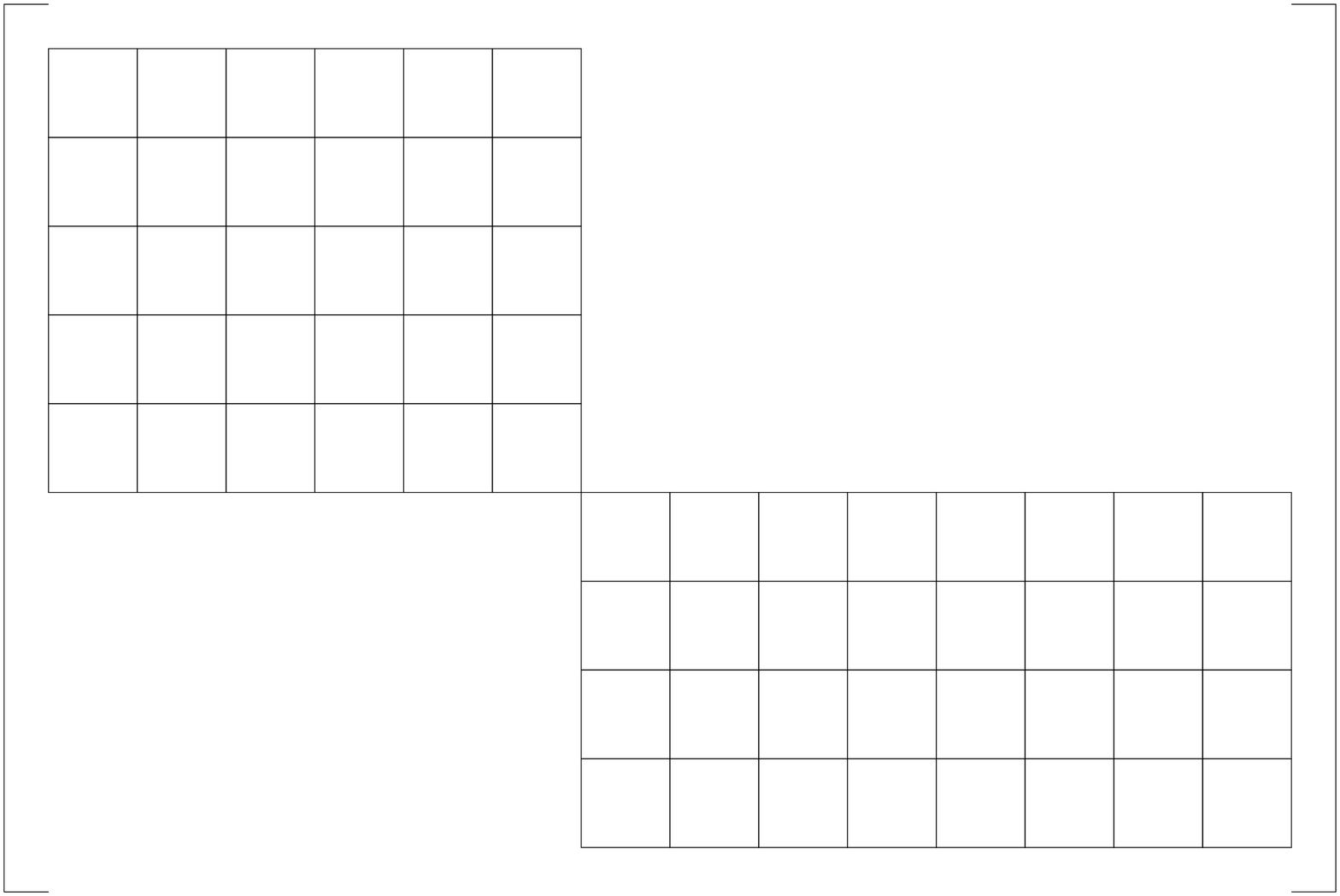}
  \raisebox{0.45in} {$+\;6x$}~\raisebox{0.05in}{\includegraphics[totalheight=0.9in]{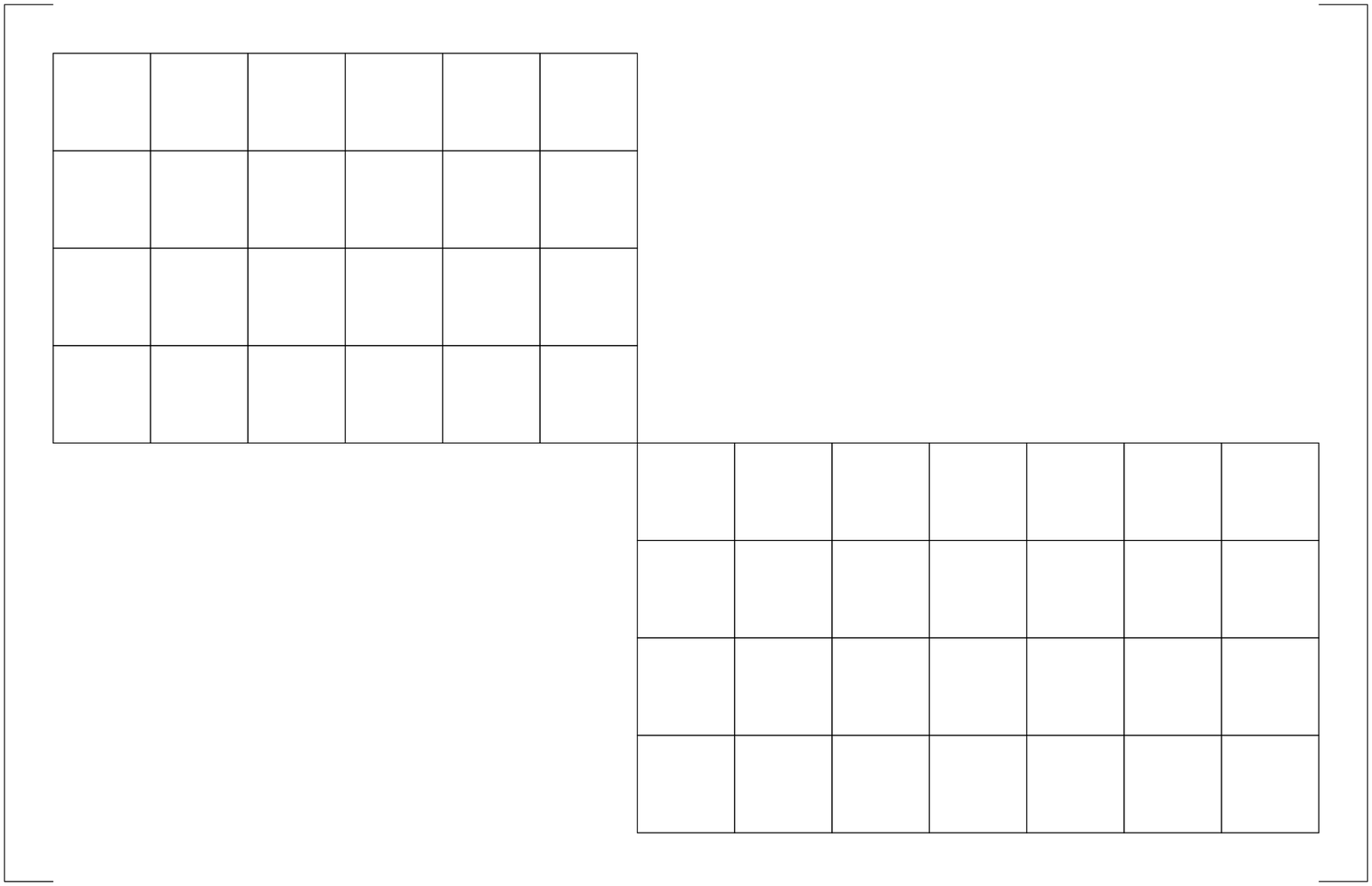}}
  \raisebox{0.45in} {$+\;6x^2$}~\raisebox{0.10in}{\includegraphics[totalheight=0.8in]{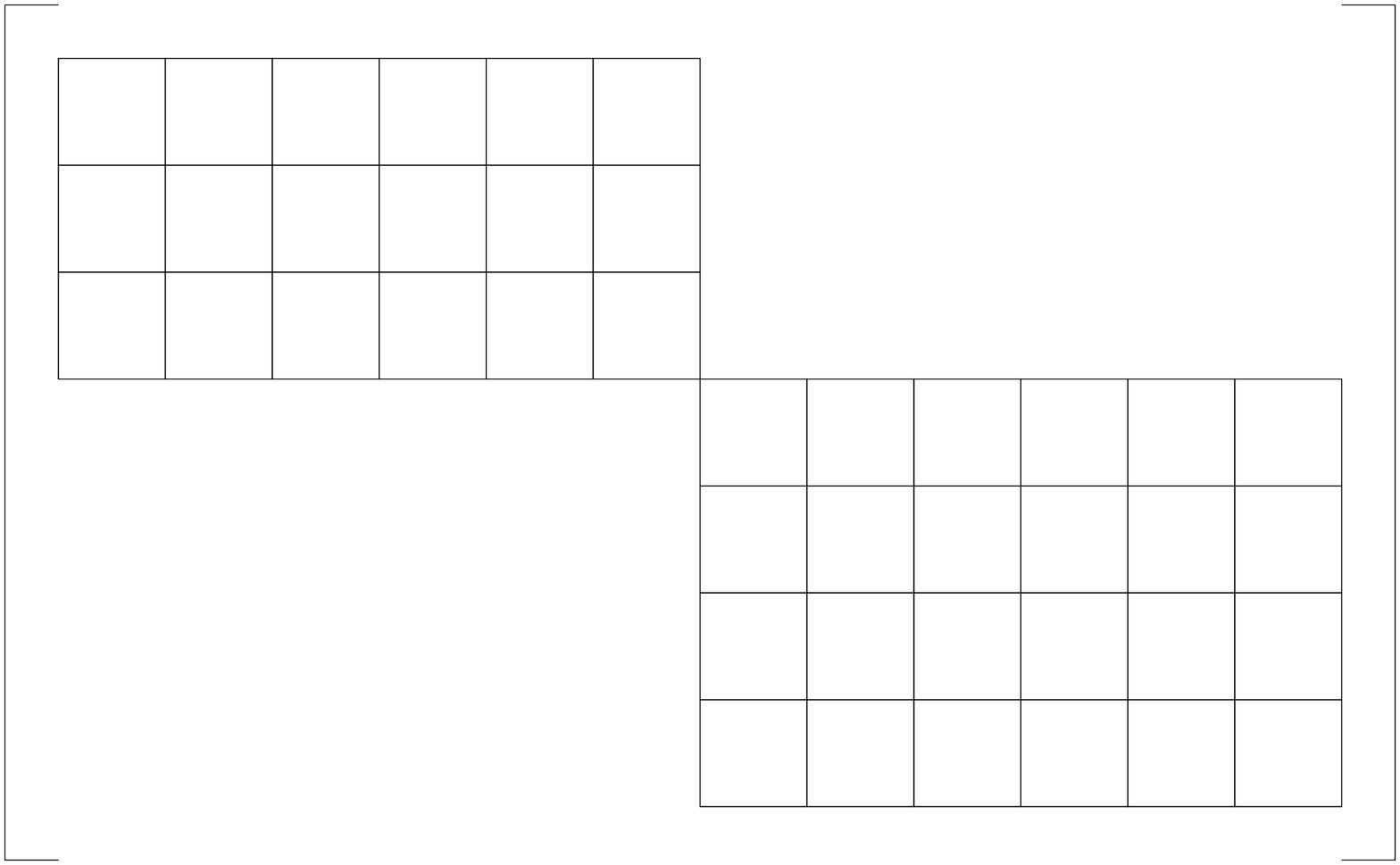}}
\end{center}

All of these are disjoint unions of rectangular boards; therefore, by Remark~\ref{dbt},
\[
[B] = R_{5,6}R_{4,8}+5xR_{4,6}R_{4,7}+5x^2R_{3,6}R_{4,6},
\]
all of which are immediately given by Theorem~\ref{rbt}.  Therefore, in this example decomposing by the largest available block was \emph{not} the most direct route to a solution.  The question of how to make optimal block choices is a complicated one, and deserves further investigation.

\end{exa}


\end{document}